\documentclass[11pt]{article}
\usepackage{amssymb,amsmath}
\usepackage{color}
\usepackage[OT2, T1]{fontenc} 
\usepackage[utf8]{inputenc} 

\usepackage{hyperref}

\newcommand{\nequiv}{\not \equiv}

\usepackage[italian, english]{babel} 
\hypersetup{hidelinks} 

\allowdisplaybreaks 

\def\claim#1.{\noindent {\bf #1.}}

\def\flushright#1{{\unskip\nobreak\hfil\penalty50\hskip2em\hbox{}\nobreak\hfil%
#1\parfillskip=0pt\finalhyphendemerits=0\par}}

\def\bull{\vrule height 1.8ex width 1.0ex depth .1ex }
\def\QED{\ifmmode\eqno\hbox{$\bull$}\else\flushright{\hbox{$\bull$}}\fi}

\newcommand{\parag}[1]{\left\{ \begin{aligned} #1 \end{aligned}\right.} 

\usepackage{scalerel}
\usepackage{tikz}
\usetikzlibrary{svg.path}
\definecolor{orcidlogocol}{HTML}{A6CE39}
\tikzset{
 orcidlogo/.pic={
 \fill[orcidlogocol] svg{M256,128c0,70.7-57.3,128-128,128C57.3,256,0,198.7,0,128C0,57.3,57.3,0,128,0C198.7,0,256,57.3,256,128z};
 \fill[white] svg{M86.3,186.2H70.9V79.1h15.4v48.4V186.2z}
 svg{M108.9,79.1h41.6c39.6,0,57,28.3,57,53.6c0,27.5-21.5,53.6-56.8,53.6h-41.8V79.1z M124.3,172.4h24.5c34.9,0,42.9-26.5,42.9-39.7c0-21.5-13.7-39.7-43.7-39.7h-23.7V172.4z}
 svg{M88.7,56.8c0,5.5-4.5,10.1-10.1,10.1c-5.6,0-10.1-4.6-10.1-10.1c0-5.6,4.5-10.1,10.1-10.1C84.2,46.7,88.7,51.3,88.7,56.8z};
 }
}
\newcommand\orcidicon[1]{\href{https://orcid.org/#1}{\mbox{\scalerel*{
\begin{tikzpicture}[yscale=-1,transform shape]
\pic{orcidlogo};
\end{tikzpicture}
}{|}}}}

\newtheorem{Theorem}{Theorem}[section]
\newtheorem{Proposition}[Theorem]{Proposition}

\newtheorem{Lemma}[Theorem]{Lemma}
\newtheorem{Remark}[Theorem]{Remark}

\newcommand{\R}{\mathbb{R}}
\newcommand{\N}{\mathbb{N}}

\newcommand{\mc}[1]{\mathcal{#1}}
\newcommand{\wto}{\rightharpoonup}

\def\half{{1\over 2}}

\def\abs#1{|#1|}

\def\norm#1{\|#1\|}

\def\epsilon{\varepsilon}
\def\eps{\varepsilon}

\def\supp{\mathop{\rm supp}}

\begin{document}

\title 
{On fractional Schrödinger equations \\ with Hartree type nonlinearities}

\author{
		Silvia Cingolani 
		\orcidicon{0000-0002-3680-9106}
		\\ \normalsize{Dipartimento di Matematica}
		\\ \normalsize{Universit\`{a} degli Studi di Bari Aldo Moro}
		\\ \normalsize{Via E. Orabona 4, 70125 Bari, Italy}
		\\ \normalsize{\href{mailto:silvia.cingolani@uniba.it}{silvia.cingolani@uniba.it}}
		\\ 
		\\Marco Gallo 
		\orcidicon{0000-0002-3141-9598}
		\\ \normalsize{Dipartimento di Matematica}
		\\ \normalsize{Universit\`{a} degli Studi di Bari Aldo Moro}
		\\ \normalsize{Via E. Orabona 4, 70125 Bari, Italy}
		\\ \normalsize{\href{mailto:marco.gallo@uniba.it}{marco.gallo@uniba.it}}
		\\ 
		\\Kazunaga Tanaka 
		\orcidicon{0000-0002-1144-1536}
		\\ \normalsize{Department of Mathematics}
		\\ \normalsize{School of Science and Engineering}
		\\ \normalsize{Waseda University}
		\\ \normalsize{3-4-1 Ohkubo, Shijuku-ku, Tokyo 169-8555, Japan}
		\\ \normalsize{\href{mailto:kazunaga@waseda.jp}{kazunaga@waseda.jp}}
	}

\date{}

\maketitle


\begin{abstract} 
\noindent Goal of this paper is to study the following doubly nonlocal equation
\begin{equation}\label{eq_abstract}
(- \Delta)^s u + \mu u = (I_\alpha*F(u))F'(u) \quad \hbox{in $\mathbb{R}^N$} \tag{P}
\end{equation}
in the case of general nonlinearities $F \in C^1(\mathbb{R})$ of Berestycki-Lions type, when $N \geq 2$ and $\mu>0$ is fixed. Here $(-\Delta)^s$, $s \in (0,1)$, denotes the fractional Laplacian, while the Hartree-type term is given by convolution with the Riesz potential $I_{\alpha}$, $\alpha \in (0,N)$. 
We prove existence of ground states of \eqref{eq_abstract}. Furthermore we obtain regularity and asymptotic decay of general solutions, extending some results contained in \cite{DSS1, MS1}.
\end{abstract}

 \medskip


\medskip



\newpage
\tableofcontents

\bigskip

\setcounter{equation}{0}
\section{Introduction}\label{section:1}

In this paper we deal with the following fractional Choquard equation
\begin{equation}\label{eq_introduction}
(- \Delta)^s u + \mu u = (I_\alpha*F(u))f(u) \quad \hbox{in $\mathbb{R}^N$} 
\end{equation}
where $N\geq 2$, $\mu>0$, $s \in (0,1)$, $\alpha \in (0,N)$, $(-\Delta)^s$ and $I_{\alpha}$ denote respectively the fractional Laplacian and the Riesz potential defined by
$$(-\Delta)^s u(x):= C_{N,s} \int_{\R^N} \frac{u(x)-u(y)}{|x-y|^{N+2s}} \, dy, \quad I_\alpha(x):= A_{N,\alpha} \frac{1}{ |x|^{N- \alpha}},$$
where $C_{N,s}:=\frac{4^s \Gamma(\frac{N+2s}{2})}{\pi^{N/2}|\Gamma(-s)|}$ and $A_{N, \alpha} := \frac{\Gamma(\frac{N-\alpha}{2})}{2^\alpha\pi^{N/2} \Gamma(\frac{\alpha}{2}) }$ are two suitable positive constants and the integral is in the principal value sense. 
Finally $F:\R \to \R$,
$F'=f$ is a nonlinearity satisfying general assumptions specified below.

When dealing with double nonlocalities, important applications arise in the study of exotic stars: minimization properties related to \eqref{eq_introduction} play indeed a fundamental role in the mathematical description of the gravitational collapse of boson stars \cite{FJL, LY} and white dwarf stars \cite{HLLS}. In fact, the study of the ground states to \eqref{eq_introduction} gives information on the size of the critical initial conditions for the solutions of the corresponding pseudo-relativistic equation \cite{Len0}. 
Moreover, when $s=\frac{1}{2}$, $N=3$, $\alpha=2$ and $F(t)=\frac{1}{r} |t|^r$, we obtain
$$\sqrt{-\Delta}u + \mu u = \left(\frac{1}{2 \pi r|x|}*|u|^r\right) |u|^{r-2} u \quad \hbox{in $\R^3$}$$
related to the well-known \emph{massless boson stars equation}~\cite{FraL,LeLe,HL}, where the pseudorelativistic operator $\sqrt{-\Delta + m}$ collapses to the square root of the Laplacian. 
Other applications can be found in relativistic physics and in quantum chemistry \cite{AM,DOS,HMT} and in the study of graphene \cite{LMM}, where the nonlocal nonlinearity describes the short time interactions between particles. 

\smallskip

In the limiting local case $s=1$, when $N=3$, $\alpha=2$ and $F(t)=\frac{1}{2} |t|^2$, 
the equation has been introduced in 1954 by Pekar in \cite{Pek0} to describe the quantum theory of a polaron at rest. Successively, in 1976 it was arisen in the work \cite{Lie1} suggested by Choquard on the modeling of an electron trapped in its own hole, in a certain approximation to Hartree-Fock theory of one-component plasma (see also \cite{FL0,FTY,Stu2}). 
In 1996 the same equation was derived by Penrose in his discussion on the self-gravitational collapse of a quantum mechanical wave-function \cite{Pen1,Pen2,Pen3,MPT} (see also \cite{TM,Tod0}) and in that context it is referred as \emph{Schr\"odinger-Newton system}.
Variational methods were also employed to derive existence and qualitative results of standing wave solutions for more generic values of $\alpha \in (0,N)$ and of power type nonlinearities $F(t)= \frac{1}{r} |t|^r$ \cite{MS0} (see also \cite{MS3, MZ,CCS1,CS0,Len1, MS1}). The case of general functions $F$, almost optimal in the sense of Berestycki-Lions \cite{BL1}, has been treated in \cite{MS1, CGT4}.

The fractional power of the Laplacian appearing in \eqref{eq_introduction}, when $s \in (0,1)$, has been introduced instead by Laskin \cite{Las0} as an extension of the classical local Laplacian in the study of nonlinear Schr\"odinger equations, replacing the path integral over Brownian motions with L\'evy flights. This operator arises naturally in many contexts and concrete applications in various fields, such as optimization, finance, crystal dislocations, charge transport in biopolymers, flame propagation, minimal surfaces, water waves, geo-hydrology, anomalous diffusion, neural systems, phase transition and Bose-Einstein condensation (see \cite{KSM,BuV,FJL,DPV,KLS,Lon0} and references therein). 
Equations involving the fractional Laplacian together with local nonlinearities have been largely investigated, and some fundamental contributions can be found in \cite{CafSil,CabSir,FLS}. In particular, existence and qualitative properties of the solutions for general classes of fractional NLS equations with local sources have been studied in \cite{FQT,CW,BKS,Iko0}.

\smallskip

Mathematically, doubly nonlocal equations have been treated in \cite{DSS1,DSS2} in the case of pure power nonlinearities (see also \cite{CFHMT} for some orbital stability results and \cite{CHHO} for a Strichartz estimates approach), obtaining existence and qualitative properties of the solutions. Other results can be found in \cite{SGY, BBMP, Luo0} for superlinear nonlinearities, in \cite{GZ} for $L^2$-supercritical Cauchy problems, in \cite{GDS} for bounded domains and in \cite{YZ} for concentration phenomena with strictly noncritical and monotone sources.

\medskip

In the present paper we address the study of \eqref{eq_introduction} when $f$ satisfies the following set of assumptions of Berestycki-Lions type \cite{BL1}:
\begin{itemize}
\item[(f1)] $f \in C(\R, \R)$;
\item[(f2)] we have
 $$i) \; \limsup_{t \to 0} \frac{|tf(t)|}{|t|^{\frac{N+ \alpha}{N}}} <+\infty, \quad
 ii) \; \limsup_{ |t| \to + \infty} \frac{|t f(t)|}{|t|^{\frac{N+ \alpha}{N-2s}}} <+\infty;$$
\item[(f3)] $F(t)= \int_0^t f(\tau) d\tau$ satisfies
 $$i) \; \lim_{t \to 0} \frac{F(t)}{|t|^{\frac{N+ \alpha}{N}}} =0, \quad
 ii) \; \lim_{ |t| \to + \infty} \frac{F(t)}{|t|^{\frac{N+ \alpha}{N-2s}}} =0;$$
\item[(f4)] there exists $t_0 \in \R$, $t_0 \neq 0$ such that $F(t_0) \neq 0$.
\end{itemize}
We observe that (f3) implies that we are in a \emph{noncritical} setting: indeed the exponents $\frac{N+\alpha}{N}$ and $\frac{N+\alpha}{N-2s}$ 
have been addressed in \cite{MS0} as critical for Choquard-type equations when $s=1$, and then generalized to $s\in (0,1)$ in \cite{DSS1}; we will assume the noncriticality in order to obtain the existence of a solution, while most of the qualitative results will be given in a \emph{possibly critical} setting. This kind of general nonlinearities include some particular cases such as pure powers $f(t) \sim t^r$, cooperating powers $f(t) \sim t^r + t^h$, competing powers $f(t)\sim t^r - t^h$ and saturable functions $f(t) \sim \frac{t^3}{1+t^2}$ (which arise, for instance, in nonlinear optics \cite{DLWZPH}).

 \smallskip

We deal first with existence of a ground state for \eqref{eq_introduction}, obtaining the following result.

\begin{Theorem}\label{th_INT_exist_unconstrained}
Assume \textnormal{(f1)--(f4)}. 
Then there exists a radially symmetric weak solution $u$ of \eqref{problem_x}, which satisfies the Pohozaev identity:
\begin{equation}\label{eq_INT_Pohozaev}
\frac{N-2s}{2}\int_{\R^N} |(-\Delta)^{s/2}u|^2 \, dx + {N \over 2} \mu \int_{\R^N} u^2 \, dx = \frac{N + \alpha}{2} \int_{\R^N} (I_{\alpha}*F(u)) F(u) \, dx .
\end{equation}
This solution is of Mountain Pass type and minimizes the energy among all the solutions satisfying \eqref{eq_INT_Pohozaev}.
\end{Theorem}

We refer to Section \ref{sec_unconstrained} for the precise meaning of \emph{weak solution}, of \emph{Mountain Pass type} and \emph{energy}, according to a variational formulation of the problem.

We point out some difficulties which arise in this framework. Indeed, the presence of the fractional power of the Laplacian does not allow to use the fact that every solution satisfies the Pohozaev identity to conclude that a Mountain Pass solution is actually a (Pohozaev) ground state, as in \cite{JT} (see Remark \ref{rem_Pohozaev_c_p}). On the other hand, the presence of the Choquard term, which scales differently from the $L^2$-norm term, does not allow to implement the classical minimization argument by \cite{CGM, BL1}. Finally, the nonhomogeneity of the nonlinearity $f$ obstructs the minimization approach of \cite{MS1,DSS1}. Thus, we need a new approach to get existence of solutions, in the spirit of \cite{CGT1, CGT2, CT1}.

Under (f1)--(f4) it is moreover possible to state the existence of a constant sign solution (see Proposition \ref{prop_ex_posiv_solut}). 
This motivates the investigation of qualitative properties for general positive solutions; in this case we consider weaker or stronger assumptions in substitution to (f1)--(f3), depending on the result. In particular, we observe that (f1)-(f2) alone imply
 $$|tf(t)| \leq C \Big(|t|^{\frac{N + \alpha}{N}} + |t|^{\frac{N + \alpha}{N-2s}}\Big),$$
 and
$$|F(t)| \leq C \Big(|t|^{\frac{N + \alpha}{N}} + |t|^{\frac{N + \alpha}{N-2s}}\Big),$$
where we notice that the last inequality is weaker than (f3); some of the qualitative results are still valid when $F$ has this possible \emph{critical} growth. 
Consider finally the following stronger assumption in the origin: 
\begin{itemize}
\item[(f5)] $\limsup_{t \to 0} \frac{|tf(t)|}{|t|^{2}} <+\infty$,
\end{itemize}
and observe that 
$$\hbox{(f5) $\implies$ (f2,i) and (f3,i).}$$

The main qualitative results that we obtain are the following ones.
\begin{Theorem}\label{th_INT_regular}
Assume \textnormal{(f1)}-\textnormal{(f2)}. Let $u\in H^s(\R^N)$ be a weak positive solution of \eqref{eq_introduction}.
Then $u \in L^1(\R^N) \cap L^{\infty}(\R^N)$. 
The same conclusion holds for generally signed solutions by assuming also \textnormal{(f5)}.
\end{Theorem}

The condition in zero of the function $f$ assumed in (f5) leads also to the following polynomial decay of the solutions.
\begin{Theorem}\label{th_INT_decay}
Assume \textnormal{(f1)}-\textnormal{(f2)} and \textnormal{(f5)}. Let $u\in H^s(\R^N)$ be a positive weak solution of \eqref{eq_introduction}. Then there exists $C', C''>0$ such that
$$\frac{C'}{1+|x|^{N+2s}} \leq u(x)\leq \frac{ C''}{1+|x|^{N+2s}},\quad \textit{ for $x \in \R^N$}.$$
\end{Theorem}

\medskip

The previous results generalize some of the ones in \cite{DSS1} to the case of general, not homogeneous, nonlinearities; in particular, we do not even assume $f$ to satisfy Ambrosetti-Rabinowitz type conditions nor monotonicity conditions. 
We observe in addition that the information $u\in L^1(\R^N)\cap L^2(\R^N)$ is new even in the power-type setting: indeed in \cite{DSS1} the authors assume the nonlinearity to be not lower critical, while here we include the possibility of criticality.
Moreover, we improve the results in \cite{SGY, Luo0} since we do not assume $f$ to be superlinear, and we have no restriction on the parameter $\alpha$.
Finally, we extend some of the results in \cite{MS1} to the fractional framework, and some of the results in \cite{BKS} to Choquard nonlinearities.

\medskip

The paper is organized as follows. We start with some notations and recalls in Section \ref{sec_prelimin}. In Section \ref{sec_unconstrained} we obtain the existence of a ground state in a noncritical setting, and in addition the existence of a positive solution. Section \ref{sec_regularity_fp} is dedicated to the study of the boundedness of positive solutions, while in Section \ref{sec_asymptotic} we investigate the asymptotic decay. Finally in the Appendix \ref{sec_bound_sign} we obtain the boundedness of general signed solutions under some more restrictive assumption.

\section{Preliminaries}\label{sec_prelimin}

Let $N\geq 2$ and $s \in (0,1)$. Recalled the definition of the fractional Laplacian \cite{DPV}
$$(-\Delta)^s u(x)= C_{N,s} \int_{\R^N} \frac{u(x)-u(y)}{|x-y|^{N+2s}} \, dy$$
for every $s \in (0,1)$, we set the fractional Sobolev space as
$$H^s(\R^N)= \big\{ u \in L^2(\R^N) \mid (-\Delta)^{s/2} u \in L^2(\R^N)\big\}$$
endowed with
$$\norm{u}_{H^s}^2= \norm{u}_2^2 + \norm{(-\Delta)^{s/2}u}_2^2.$$
In particular, we consider the subspace of radially symmetric functions $H^s_r(\R^N)$, and recall the continuous embedding \cite[Theorem 3.5]{DPV}
$$H^s(\R^N) \hookrightarrow L^p(\R^N)$$
for every $p \in [2, 2^*_s]$, $2^*_s = \frac{2N}{N-2s}$ critical Sobolev exponent, and the compact embedding \cite{Lio1}
$$H^s_r(\R^N) \hookrightarrow \hookrightarrow L^p(\R^N)$$
for every $p \in (2, 2^*_s)$. In addition we have the following embedding of the homogeneous space \cite[Theorem 6.5]{DPV} for some $\mc{S}>0$
\begin{equation}\label{eq_embd_homog}
\norm{u}_{2^*_s} \leq \mc{S}^{-1/2} \norm{(-\Delta)^{s}u}_2.
\end{equation}

Moreover the following relation with the Gagliardo seminorm holds \cite[Proposition 3.6]{DPV}, for some $C(N,s)>0$
\begin{equation}\label{eq_semin_gagl}
 \norm{(-\Delta)^{s/2} u}_2^2 = C(N,S) \int_{\R^{2N}} \frac{|u(x)-u(y)|^2}{|x-y|^{N+2s}} \, dx \, dy.
\end{equation}
Thanks to this last formulation, we obtain that if $u\in H^s(\R^N)$ and $h: \R \to \R$ is a Lipschitz function with $h(0)=0$, then $h(u) \in H^s(\R^N)$. Indeed
$$\norm{h(u)}_2^2 = \int_{\R^N}|h(u)-h(0)|^2 \, dx \leq \int_{\R^N} \norm{h'}_{\infty}^2 |u-0|^2 \, dx = \norm{h'}_{\infty}^2 \norm{u}_2^2$$
and
$$\norm{(-\Delta)^{s/2} h(u)}_2^2 \leq C(N,S) \int_{\R^{2N}} \frac{\norm{h'}_{\infty}^2 |u(x)-u(y)|^2}{|x-y|^{N+2s}} \, dx \, dy = \norm{h'}_{\infty}^2 \norm{(-\Delta)^{s/2} u}_2^2.$$

We further have the following relation with the Fourier transform \cite[Proposition 3.3]{DPV}
$$(-\Delta)^s u = \mc{F}^{-1}(|\xi|^{2s}(\mc{F}(u));$$
notice that this last expression is suitable for defining the fractional Sobolev space $W^{s,p}(\R^N)$ also for $s\geq 1$ and $p\geq 1$, by \cite{FQT}
$$W^{s,p}(\R^N)= \big\{ u \in L^p(\R^N) \mid \mc{F}^{-1}(|\xi|^{s}(\mc{F}(u)) \in L^p(\R^N)\big\}.$$
 
Finally, set $\alpha \in (0,N)$, we recall the following standard estimates for the Riesz potential \cite[Theorem 4.3]{LiLo}.

\begin{Proposition}[Hardy-Littlewood-Sobolev inequality]\label{prop_HLS}
Let $\alpha \in (0, N)$, and let $r, \,h \in (1,+\infty)$ be such that $\frac{1}{r} - \frac{1}{h} = \frac{\alpha}{N}$. Then the map
$$f \in L^r(\R^N) \mapsto I_{\alpha}* f \in L^h(\R^N) $$
is continuous. 
In particular, if $r, \, t \in (1, +\infty)$ verify $\frac 1 r + \frac 1 t = \frac{N+\alpha}N$, then there exists a constant $C=C(N,\alpha,r,t)>0$ such that 
$$
\left| \int_{\R^N} (I_\alpha*g) h\, dx \right| \leq C\|g\|_r \|h\|_{t}
$$
for all $g\in L^r(\R^N)$ and $h\in L^t(\R^N)$. 
\end{Proposition}

\section{Existence of ground states} \label{sec_unconstrained}

In this section we search for solutions to the fractional Choquard equation
\begin{equation} \label{problem_x}
(- \Delta)^s u + \mu u =	(I_\alpha*F(u))f(u) \quad \hbox{in $\mathbb{R}^N$}
\end{equation}
by variational methods on the subspace of radially symmetric functions $H^s_r(\R^N)$. We recall that $F'=f$ and we assume (f1)-(f2) in order to have well defined functionals. 
We set $\mc{D}: H^s_r(\R^N) \to \R$ as
$$ \mc{D}(u) := \int_{\R^N} (I_\alpha*F(u))F(u)\, dx$$
and define the $C^1$-functional $\mc{J}_{\mu}: H^s_r(\R^N) \to \R$ associated to \eqref{problem_x} by 
$$\mc{J}_{\mu}(u):= \half \int_{\R^N} |(-\Delta)^{s/2} u|^2\,dx -\half {\mathcal D}(u) + \frac{\mu}{2} \|u\|_2^2.$$
We notice that, by the Principle of Symmetric Criticality of Palais, the critical points of $\mc{J}_{\mu}$ are \emph{weak solutions} of \eqref{problem_x}.
Moreover, inspired by the Pohozaev identity
\begin{equation}\label{eq_Poh_ident}
\frac{N-2s}{2}\norm{(-\Delta)^{s/2}u}_{2}^2 + {N \over 2} \mu \|u\|_2^2 = \frac{N + \alpha}{2} \, {\mathcal D}(u)
\end{equation}
we define also the Pohozaev functional $\mc{P}_{\mu}: H^s_r(\R^N) \to \R$ by 
$$\mc{P}_{\mu}(u):=\frac{N-2s}{2}\norm{(-\Delta)^{s/2}u}_{2}^2 -\frac{N+ \alpha}{2}{\mathcal D}(u)+ \frac{N}{2} \mu \|u\|_2^2.$$
Furthermore we introduce the set of paths
$$\Gamma_{\mu} := \big\{ \gamma \in C\big([0,1], H^s_r(\R^N)\big) \mid \gamma(0)=0, \, \mc{J}_{\mu}(\gamma(1))<0\big\}$$
and the \emph{Mountain Pass} (MP for short) value 
\begin{equation}\label{eq_MP_value}
 l (\mu) := \inf_{\gamma\in\Gamma_\mu}\max_{t\in [0,1]} \mc{J}_{\mu}(\gamma(t)).
 \end{equation}
Finally we set
$$p(\mu):= \inf \big\{ \mc{J}_{\mu}(u) \mid u \in H^s_r(\R^N) \setminus \{0\}, \; \mc{P}_{\mu}(u)=0\big\}$$
the \emph{least energy} of $\mc{J}_{\mu}$ on the Pohozaev set.

\begin{Remark}\label{rem_buona_posit}
Since of key importance in the good definition of the functionals, as well as in bootstrap argument in the rest of the paper, we write here in which spaces lie the considered quantities. 
Let $u\in H^s(\R^N) \subset L^2(\R^N) \cap L^{2^*_s}(\R^N)$. By \textnormal{(f2)} we have
\begin{align*}
f(u) &\in L^{\frac{2N}{\alpha}}(\R^N)\cap L^{\frac{N}{\alpha} \frac{2N}{N-2s}}(\R^N)+ L^{2 \frac{N-2s}{\alpha+2s}}\cap L^{\frac{2N}{\alpha+2s}}(\R^N) \\
&\subset L^{\frac{2N}{\alpha}}(\R^N) + L^{\frac{2N}{\alpha+2s}}(\R^N),\\
F(u) &\in L^{\frac{2N}{N+\alpha}}(\R^N)\cap L^{\frac{N}{N+\alpha} \frac{2N}{N-2s}}(\R^N) + L^{2 \frac{N-2s}{N+\alpha}}(\R^N)\cap L^{\frac{2N}{N+\alpha}}(\R^N) \\
&\subset L^{\frac{2N}{N+\alpha}}(\R^N).
\end{align*}
Thus by the Hardy-Littlewood-Sobolev inequality we obtain
\begin{align*}
I_{\alpha}*F(u) &\in L^{\frac{2N}{N-\alpha}}(\R^N) \cap L^{\frac{2N^2}{N^2-(\alpha+2s)N-2s\alpha}}(\R^N) + L^{\frac{2N(N-2s)}{N^2-\alpha N+4s \alpha}}(\R^N) \cap L^{\frac{2N}{N-\alpha}}(\R^N) \\
&\subset L^{\frac{2N}{N-\alpha}}(\R^N).
\end{align*}
Finally, by the H\"older inequality, we have
\begin{align*}
(I_{\alpha}*F(u))f(u) &\in L^2(\R^N)\cap L^{\frac{2N^2}{N^2 - 2s\alpha}}(\R^N) + L^{\frac{2N(N-2s)}{N^2+2 \alpha s}}(\R^N)\cap L^{\frac{2N}{N+2s}}(\R^N)\\
& \subset L^2(\R^N) + L^{\frac{2N}{N+2s}}(\R^N).
\end{align*}
In particular we observe that $(I_{\alpha}*F(u))f(u)$ does not lie in $L^2(\R^N)$, generally. On the other hand, if $\varphi \in H^s(\R^N) \subset L^2(\R^N) \cap L^{2^*_s}(\R^N)$, we notice that the found summability of $(I_{\alpha}*F(u))f(u)$ is enough to have
$$\int_{\R^N} (I_{\alpha}*F(u))f(u) \varphi \, dx$$
well defined. 
\end{Remark}

We present now an existence result for \eqref{problem_x}.

\begin{Theorem}\label{th_exist_unconstrained}
Assume \textnormal{(f1)--(f4)}. Let $\mu>0$ be fixed. Then there exists a Mountain Pass solution $u$ of \eqref{problem_x}, that is
$$\mc{J}_{\mu}(u)=l(\mu)>0.$$
Moreover, the found solution satisfies the Pohozaev identity
$$\mc{P}_{\mu}(u)=0.$$
\end{Theorem}

\claim Proof.
We split the proof in some steps.

\textbf{Step 1.}
We first show that $\mc{J}_{\mu}$ satisfies the Palais-Smale-Pohozaev condition at every level $b \in \R$, that is each sequence $u_n$ in $H^s_r(\R^N)$ satisfying
\begin{equation}\label{eq_dim_u_1}
\mc{J}_{\mu}(u_n) \to b,
\end{equation}
\begin{equation}\label{eq_dim_u_2}
\mc{J}'_{\mu}( u_n) \to 0 \quad \ \hbox{strongly in $(H^s_r(\R^N))^*$},
\end{equation}
\begin{equation}\label{eq_dim_u_3}
\mc{P}_{\mu}(u_n) \to 0,
\end{equation}
converges up to a subsequence. Indeed \eqref{eq_dim_u_1} and \eqref{eq_dim_u_3} imply
$$\frac{\alpha +2s}{2} \norm{(-\Delta)^{s/2} u_n}_2^2 + \frac{\alpha}{2} \mu \norm{u_n}_2^2 = (N+\alpha)b + o(1).$$
Thus we obtain that $b\geq 0$ and $u_n$ is bounded in $H^s_r(\R^N)$. 

\textbf{Step 2.} 
After extracting a subsequence, denoted in the same way, we may assume that $u_n \wto u_0$ weakly in $H^s_r(\R^N)$.
Taking into account the assumptions (f1)--(f3), we obtain
$$
\int_{\R^N} (I_\alpha \ast F(u_n)) f(u_n) u_0 \, dx
\to 
\int_{\R^N} (I_\alpha \ast F(u_0)) f(u_0) u_0 \, dx
$$
and 
$$
\int_{\R^N} (I_\alpha \ast F(u_n)) f(u_n) u_n\, dx
\to 
\int_{\R^N} (I_\alpha \ast F(u_0)) f(u_0) u_0\, dx.
$$
Thus we derive that 
$
\langle \mc{J}_{\mu}'(u_n), u_n \rangle \to 0$ 
and 
$
\langle \mc{J}_{\mu}'( u_n), u_0 \rangle \to 0$,
and hence
$$\|(- \Delta)^{s/2} u_n \|_2^2 +\mu \|u_n\|^2_2 \to 
\|(- \Delta)^{s/2} u_0 \|^2_2
 + \mu \|u_0\|^2_2 
$$ 
which implies $u_n \to u_0$ strongly in $H^s_r(\R^N)$.

\textbf{Step 3.} Denote by 
$$[\mc{J}_{\mu} \leq b]:= \big\{u\in H^s_r(\R^N) \mid \mc{J}_{\mu}(u) \leq b \big\}$$
the sublevel of $\mc{J}_{\mu}$ and by 
$$K_b:= \big\{ u \in H^s_r(\R^N)\mid \mc{J}_{\mu}(u)=b, \, \mc{J}_{\mu}'(u)=0, \, \mc{P}_{\mu}(u)=0 \big\}$$
 the set of critical points of $\mc{J}_{\mu}$ satisfying the Pohozaev identity. Then, by Steps 1-2, $K_b$ is compact. 
Arguing as in \cite[Proposition 4.5]{HT} (see also \cite[Proposition 3.1 and Corollary 4.3]{IT0}), we obtain for any $b\in \R$, $\bar{\eps}>0$ and any $U$ open neighborhood of $K_b$, that there exist an $\epsilon \in (0, \bar \epsilon)$ and a continuous map $\eta:[0,1] \times H^s_r(\R^N) \to H^s_r(\R^N)$ such that 
\begin{itemize}
	\item[$(1^o)$] 
	$
	 \eta(0,u)=u \quad \forall u \in H^s_r(\R^N)
	$;
	\item[$(2^o)$] 
		$
	\eta(t, u)=u \quad \forall (t,u) \in [0,1] \times [\mc{J}_{\mu}\leq b- \bar\epsilon]
	$;
	\item[$(3^o)$] 
	$
	\mc{J}_{\mu}( \eta(t, u)) \leq \mc{J}_{\mu}(u) \quad 
	\ \forall (t,u) \in [0,1] \times H^s_r(\R^N)
	$;
	\item[$(4^o)$] 
	$ \eta (1, [\mc{J}_{\mu} \leq b+ \epsilon] \setminus U) \subset [\mc{J}_{\mu} \leq b- \epsilon]
	$;
				\item[$(5^o)$] 
	$ \eta (1, [\mc{J}_{\mu} \leq b+ \epsilon]) \subset [\mc{J}_{\mu} \leq b - \epsilon] \cup U
	$;
	\item[$(6^o)$] 
	if $K_b= \emptyset$, then
	$ \eta (1, [\mc{J}_{\mu} \leq b+ \epsilon]) \subset [\mc{J}_{\mu} \leq b- \epsilon].
	$
\end{itemize}

\textbf{Step 4.} By exploiting (f4) and arguing as in \cite[Proposition 2.1]{MS1}, we obtain the existence of a function $v \in H^s_r(\R^N)$ such that $\mc{D}(v)>0$. Thus defined $\gamma(t) := v(\cdot/t)$ for $t > 0$ and $\gamma(0): = 0$ we have $\mc{J}(\gamma(t))<0$ for $t$ large and $\mc{J}(\gamma(t))>0$ for $t$ small; this means, after a suitable rescaling, that $l(\mu)$ is finite and strictly positive. 
In particular we observe that $0\notin K_{l(\mu)}$.

\textbf{Step 5.} By applying the deformation result at level $b= l(\mu)>0$, the existence of a Mountain Pass solution $u$ is then obtained classically. 
Moreover, $u \in K_{l(\mu)}$ by construction, thus $u\nequiv 0$ and $\mc{P}_{\mu}(u)=0$.
\QED

\bigskip

We prove now that the found solution is actually a ground state over the Pohozaev set.

\begin{Proposition}\label{prop_MP=PM}
The Mountain Pass level and the Pohozaev minimum level coincide, that is
$$l(\mu)=p(\mu)>0.$$
In particular, the solution found in Theorem \ref{th_exist_unconstrained} is a Pohozaev minimum.
\end{Proposition}

\claim Proof.
Let $u \in H^s_r(\R^N)\setminus\{0\}$ such that $\mc{P}_{\mu}(u)=0$; observe that $\mc{D}(u)>0$. We define $\gamma(t):=u(\cdot/t)$ for $t \neq 0$ and $\gamma(0):=0$ so that $t\in (0,+\infty) \mapsto \mc{J}_{\mu}(\gamma(t))$ is negative for large values of $t$, and it attains the maximum in $t=1$. After a suitable rescaling we have $\gamma \in \Gamma_{\mu}$ and thus
\begin{equation}\label{eq_dis_P-MP}
\mc{J}_{\mu}(u) = \max_{t \in [0,1]}\mc{J}_{\mu}(\gamma(t)) \geq l(\mu).
\end{equation}
Passing to the infimum in equation \eqref{eq_dis_P-MP} we have $p(\mu) \geq l(\mu)$. Let now $\gamma \in \Gamma_{\mu}$. By definition we have $\mc{J}_{\mu}(\gamma(1))<0$, thus by
$$ \mc{P}_{\mu}(v) = N \mc{J}_{\mu}(v) - s\norm{(-\Delta)^{s/2} v}_2^2 - \frac{\alpha}{2} \mc{D}(v), \quad v \in H^s_r(\R^N),$$
we obtain $\mc{P}_{\mu}(\gamma(1))<0$. 
In addition, since $\mc{D}(u)=o(\norm u_{H^s}^2)$ as $u \to 0$ and $\gamma(t)\to 0$ as $t \to 0$ in $H^s_r(\R^N)$, we have
$$\mc{P}_{\mu}(\gamma(t))>0 \quad \hbox{for small $t>0$}.$$
Thus there exists a $t^*$ such that $\mc{P}_{\mu}(\gamma(t^*))=0$, and hence
$$p(\mu) \leq \mc{J}_{\mu}(\gamma(t^*)) \leq \max_{t \in [0,1]} \mc{J}_{\mu}(\gamma(t));$$
passing to the infimum we come up with $p(\mu) \leq l(\mu)$, and hence the claim.
\QED

\bigskip

\medskip

\claim Proof of Theorem \ref{th_INT_exist_unconstrained}.
We obtain the result by matching Theorem \ref{th_exist_unconstrained} and Proposition \ref{prop_MP=PM}.
\QED

\bigskip

We pass to investigate more in details Pohozaev minima, showing that it is a general fact that they are solutions of the equation \eqref{problem_x}.

\begin{Proposition}\label{prop_poho_min_sol}
Every Pohozaev minimum is a solution of \eqref{problem_x}, i.e.
$$\mc{J}_{\mu}(u)=p(\mu) \; \hbox{ and } \; \mc{P}_{\mu}(u)=0$$
imply
$$\mc{J}'_{\mu}(u)=0.$$
As a consequence
$$p(\mu)= \inf \big\{ \mc{J}_{\mu}(u) \mid u \in H^s_r(\R^N) \setminus \{0\}, \; \mc{P}_{\mu}(u)=0, \; \mc{J}'_{\mu}(u)=0\big\}.$$
\end{Proposition}

\claim Proof.
Let $u$ be such that $\mc{J}_{\mu}(u)=p(\mu)$ and $\mc{P}_{\mu}(u)=0$. 
In particular, considered $\gamma(t)=u(\cdot/t)$, we have that $\mc{J}_{\mu}(\gamma(t))$ is negative for large values of $t$ and its maximum value is $p(\mu)$ attained only in $t=1$. 

Assume by contradiction that $u$ is not critical. Let $I:=[1-\delta, 1+\delta]$ be such that $\gamma(I) \cap K_{p(\mu)}=\emptyset$, and set $\bar{\eps} := p(\mu) - \max_{t \notin I} \mc{J}_{\mu}(\gamma(t)) >0$.
Let now $U$ be a neighborhood of $K_{p(\mu)}$ verifying $\gamma(I) \cap U = \emptyset$: by the deformation lemma presented in the proof of Theorem \ref{th_exist_unconstrained} there exists an $\eta:[0,1]\times H^s_r(\R^N) \to H^s_r(\R^N)$ at level $p(\mu)\in\R$ with properties $(1^o)$-$(6^o)$.
 Define then $\tilde{\gamma}(t):= \eta(1, \gamma(t))$ a deformed path. 
 
 For $t \notin I$ we have $\mc{J}_{\mu}(\gamma(t))< p(\mu) - \bar{\eps}$, and thus by $(2^o)$ we gain
\begin{equation}\label{eq_dim_P_sol_1}
\mc{J}_{\mu}(\tilde{\gamma}(t)) = \mc{J}_{\mu}(\gamma(t)) < p(\mu) - \bar{\eps}, \quad \hbox{ for $t \notin I$}.
\end{equation}
Let now $t \in I$: we have $\gamma(t) \notin U$ and $\mc{J}_{\mu}(\gamma(t))\leq p(\mu) \leq p(\mu)+ \eps$, thus by $(4^o)$ we obtain
\begin{equation}\label{eq_dim_P_sol_2}
\mc{J}_{\mu}(\tilde{\gamma}(t)) \leq p(\mu) - \eps.
\end{equation}
 Joining \eqref{eq_dim_P_sol_1} and \eqref{eq_dim_P_sol_2} we have
$$\max_{t \geq 0} \mc{J}_{\mu}(\tilde{\gamma}(t)) < p(\mu)=l(\mu)$$
which is an absurd, since after a suitable rescaling it results that $\tilde{\gamma} \in \Gamma_{\mu}$, thanks to $(3^o)$. 
\QED

\medskip

\begin{Remark}\label{rem_Pohozaev_c_p}
We point out that it is not known, even in the case of local nonlinearities \cite{BKS}, if 
$$p(\mu)= \inf \big\{ \mc{J}_{\mu}(u) \mid u \in H^s_r(\R^N) \setminus \{0\}, \; \mc{J}_{\mu}'(u)=0\big\}.$$
On the other hand, by assuming that every solution of \eqref{problem_x} satisfies the Pohozaev identity (see e.g. \cite[Proposition 2]{SGY} and \cite[Equation (6.1)]{DSS1}), the claim holds true.
\end{Remark}

We show now that, under the same assumptions of Theorem \ref{th_exist_unconstrained}, we can find a solution with constant sign.

\begin{Proposition}\label{prop_ex_posiv_solut}
Assume \textnormal{(f1)--(f4)} and that $F\nequiv 0$ on $(0,+\infty)$ (i.e., $t_0$ in assumption \textnormal{(f4)} can be chosen positive). Then there exists a positive radially symmetric solution of \eqref{problem_x}, which is minimum over all the positive functions on the Pohozaev set.
\end{Proposition}

\claim Proof.
Let us define
$$g:= \chi_{(0, +\infty)} f.$$
We have that $g$ still satisfies (f1)--(f4). Thus, by Theorem \ref{th_exist_unconstrained} there exists a solution $u$ of
$$(- \Delta)^s u + \mu u = (I_\alpha*G(u))g(u) \quad \hbox{in $\mathbb{R}^N$}$$
where $G(t):=\int_0^t g(\tau) d \tau$. 
We show now that $u$ is positive. We start observing the following: by \eqref{eq_semin_gagl} we have
\begin{eqnarray*}
\norm{(-\Delta)^{s/2} |u|}_2^2 &=& C(N,s) \int_{\R^{2N}} \frac{\big(|u(x)|-|u(y)|\big)^2}{|x-y|^{N+2s}} \, dx \, dy \\
&=& C(N,s) \int_{\R^{2N}} \frac{|u|^2(x) + |u|^2(y) - 2|u|(x)|u|(y)}{|x-y|^{N+2s}} \, dx \, dy \\
&\leq& C(N,s) \int_{\R^{2N}} \frac{u^2(x) + u^2(y) - 2u(x)u(y)}{|x-y|^{N+2s}} \, dx \, dy \\
&=& C(N,s) \int_{\R^{2N}} \frac{\big(u(x)-u(y)\big)^2}{|x-y|^{N+2s}} \, dx \, dy = \norm{(-\Delta)^{s/2}u}_2^2,
\end{eqnarray*}
thus
$$\norm{(-\Delta)^{s/2} |u|}_2 \leq \norm{(-\Delta)^{s/2} u}_2.$$
In particular, written $u=u_+ - u_-$, by the previous argument we have $u_- = \frac{|u|-u}{2} \in H^s_r(\R^N)$. Thus, chosen $u_-$ as test function, we obtain
$$\int_{\R^N} (-\Delta)^{s/2} u \, (-\Delta)^{s/2}u_- \, dx + \mu \int_{\R^N} u \, u_- \, dx = \int_{\R^N} (I_{\alpha} * G(u))g(u) u_- \, dx.$$
By definition of $g$ and \eqref{eq_semin_gagl} we have
\begin{equation}\label{eq_splitting_posit}
C_{N,s}\int_{\R^N \times \R^N} \frac{(u(x)-u(y))(u_-(x)-u_-(y))}{|x-y|^{N+2s}} \, dx \, dy - \mu \int_{\R^N} u_-^2 \, dx =0.
\end{equation}
Splitting the domain, we gain
\begin{eqnarray*}
\lefteqn{\int_{\R^N \times \R^N} \frac{(u(x)-u(y))(u_-(x)-u_-(y))}{|x-y|^{N+2s}} \, dx \, dy =}\\
&&- \int_{\{u(x)\geq 0\} \times \{u(y)<0\}} \frac{(u_+(x)+u_-(y))(u_-(y))}{|x-y|^{N+2s}} \, dx \, dy -\\
&& -\int_{\{u(x)<0\} \times \{u(y)\geq 0\}} \frac{(u_-(x)+u_+(y))(u_-(x))}{|x-y|^{N+2s}} \, dx \, dy -\\
&&- \int_{\{u(x)<0\} \times \{u(y)<0\}} \frac{(u_-(x)-u_-(y))^2}{|x-y|^{N+2s}} \, dx \, dy.
\end{eqnarray*}
Thus we obtain that the left-hand side of \eqref{eq_splitting_posit} is sum of non positive pieces, thus $u_- \equiv 0$, that is $u\geq 0$. Hence $g(u)=f(u)$ and $G(u)=F(u)$, which imply that $u$ is a (positive) solution of \eqref{problem_x}.
\QED

\medskip

\section{Regularity}\label{sec_regularity_fp}

In this section we prove some regularity results for \eqref{problem_x}. 
We split the proof of Theorem \ref{th_INT_regular} in different steps.

We start from the following lemma, that can be found in \cite[Lemma 3.3]{MS1}.
\begin{Lemma}[\cite{MS1}]\label{lem_primoMS}
Let $N \geq 2$ and $\alpha \in (0, N)$. Let $\lambda \in [0,2]$ and $q,r, h,k \in [1, +\infty)$ be such that
$$1+ \frac{\alpha}{N} - \frac{1}{h} - \frac{1}{k} = \frac{\lambda}{q} + \frac{2-\lambda}{r}.$$
Let $\theta \in (0,2)$ satisfying
$$\min\{q,r\} \left(\frac{\alpha}{N} - \frac{1}{h}\right) < \theta < \max\{q,r\} \left( 1- \frac{1}{h}\right),$$
$$\min\{q,r\} \left(\frac{\alpha}{N} - \frac{1}{k}\right) < 2-\theta < \max\{q,r\} \left( 1- \frac{1}{k}\right).$$
Let $H \in L^h(\R^N)$, $K \in L^k(\R^N)$ and $u \in L^q(\R^N) \cap L^r(\R^N)$. Then
$$\int_{\R^N} \left( I_{\alpha} * \big(H|u|^{\theta}\big) \right) K|u|^{2-\theta} \, dx \leq C \norm{H}_h \norm{K}_k \norm{u}_q^{\lambda} \norm{u}_r^{2-\lambda}$$
for some $C>0$ (depending on $\theta$). 
\end{Lemma}

By a proper use of Lemma \ref{lem_primoMS} we obtain now an estimate on the Choquard term depending on $H^s$-norm of the function.

\begin{Lemma}\label{lem_secondMS}
Let $N\geq 2$, $s \in (0,1)$ and $\alpha \in (0,N)$. 
Let moreover $\theta \in (\frac{\alpha}{N}, 2- \frac{\alpha}{N})$ and $H, K \in L^{\frac{2N}{\alpha}}(\R^N) + L^{\frac{2N}{\alpha+2s}}(\R^N) $. Then for every $\eps>0$ there exists $C_{\eps, \theta}>0$ such that
$$\int_{\R^N} \left(I_{\alpha}*\big(H|u|^{\theta}\big)\right)K|u|^{2-\theta} \, dx \leq \eps^2 \norm{(-\Delta)^{s/2} u}_2^2 + C_{\eps, \theta}\norm{u}_2^2$$
for every $u \in H^s(\R^N)$.
\end{Lemma}

\claim Proof.
Observe that $2-\theta \in (\frac{\alpha}{N}, 2- \frac{\alpha}{N})$ as well. We write
$$H= H^* + H_* \in L^{\frac{2N}{\alpha}}(\R^N) + L^{\frac{2N}{\alpha+2s}}(\R^N),$$
$$K= K^* + K_* \in L^{\frac{2N}{\alpha}}(\R^N) + L^{\frac{2N}{\alpha+2s}}(\R^N).$$
We split $\int_{\R^N} \left(I_{\alpha}*\big(H|u|^{\theta}\big)\right)K|u|^{2-\theta} \, dx $ in four pieces and choose
$$q=r=2, \quad h=k=\frac{2N}{\alpha}, \quad \lambda=2,$$
$$q = 2, \; r=\frac{2N}{N-2s}, \quad h=\frac{2N}{\alpha}, \; k=\frac{2N}{\alpha+2s}, \quad \lambda=1,$$
$$q = 2, \; r=\frac{2N}{N-2s}, \quad h=\frac{2N}{\alpha+2s}, \; k=\frac{2N}{\alpha}, \quad \lambda=1,$$
$$q=r=\frac{2N}{N-2s}, \quad h=k=\frac{2N}{\alpha+2s}, \quad \lambda=0,$$
in Lemma \ref{lem_primoMS}, to obtain
\begin{align*}
\int_{\R^N} \left(I_{\alpha}*\big(H|u|^{\theta}\big)\right)K|u|^{2-\theta} \, dx \lesssim& \norm{H^*}_{\frac{2N}{\alpha}} \norm{K^*}_{\frac{2N}{\alpha}} \norm{u}_2^2 + \norm{H^*}_{\frac{2N}{\alpha}}\norm{K_*}_{\frac{2N}{\alpha+2s}} \norm{u}_2 \norm{u}_{\frac{2N}{N-2s}}+ \\
&+ \norm{H_*}_{\frac{2N}{\alpha+2s}} \norm{K^*}_{\frac{2N}{\alpha}}\norm{u}_2 \norm{u}_{\frac{2N}{N-2s}} + \norm{H_*}_{\frac{2N}{\alpha+2s}}\norm{K_*}_{\frac{2N}{\alpha+2s}} \norm{u}_{\frac{2N}{N-2s}}^2.
\end{align*}
Recalled that $\frac{2N}{N-2s}=2^*_s$ and the Sobolev embedding \eqref{eq_embd_homog}, we obtain
\begin{align}
\int_{\R^N} \left(I_{\alpha}*\big(H|u|^{\theta}\big)\right)K|u|^{2-\theta} \, dx \lesssim & \left(\norm{H^*}_{\frac{2N}{\alpha}}\norm{K^*}_{\frac{2N}{\alpha}}\right) \norm{u}_2^2 + \left(\norm{H_*}_{\frac{2N}{\alpha+2s}} \norm{K_*}_{\frac{2N}{\alpha+2s}}\right) \norm{(-\Delta)^{s/2} u}_2^2 + \notag \\
& + \left(\norm{H^*}_{\frac{2N}{\alpha}}\norm{K_*}_{\frac{2N}{\alpha+2s}} + \norm{H_*}_{\frac{2N}{\alpha+2s}} \norm{K^*}_{\frac{2N}{\alpha}}\right) \norm{u}_2 \norm{(-\Delta)^{s/2} u}_2,\label{eq_dim_mixed_terms_HK}
\end{align}
where $\lesssim$ denotes an inequality up to a constant.
We want to show now that, since $ \frac{2N}{\alpha}>\frac{2N}{\alpha+2s}$, we can choose the decomposition of $H$ and $K$ such that the $L^{\frac{2N}{\alpha+2s}}$-pieces are arbitrary small (see \cite[Lemma 2.1]{BK0}). Indeed, let
$$H = H_1 + H_2 \in L^{\frac{2N}{\alpha}}(\R^N) + L^{\frac{2N}{\alpha+2s}}(\R^N)$$
a first decomposition. Let $M>0$ to be fixed, and write
$$H = \left( H_1 + H_2 \chi_{\{|H_2| \leq M\}}\right) + H_2 \chi_{\{|H_2| >M\}}.$$
Since $ H_2 \chi_{\{|H_2| \leq M\}} \in L^{\frac{2N}{\alpha+2s}}(\R^N) \cap L^{\infty}(\R^N)$ and $\frac{2N}{\alpha} \in (\frac{2N}{\alpha+2s}, \infty)$, we have $H_2 \chi_{\{|H_2| \leq M\}} \in L^{\frac{2N}{\alpha}}(\R^N)$, and thus
$$H^*:= H_1 + H_2 \chi_{\{|H_2| \leq M\}} \in L^{\frac{2N}{\alpha}}(\R^N), \quad H_* : =H_2 \chi_{\{|H_2| >M\}} \in L^{\frac{2N}{\alpha+2s}}(\R^N).$$
On the other hand
$$\norm{H_*}_{\frac{2N}{\alpha+2s}} = \left(\int_{|H_2|>M} |H_2|^{\frac{2N}{\alpha+2s}} \, dx\right)^{\frac{\alpha+2s}{2N}}$$
which can be made arbitrary small for $M\gg 0$. In particular we choose the decomposition so that
$$\left(\norm{H_*}_{\frac{2N}{\alpha+2s}} \norm{K_*}_{\frac{2N}{\alpha+2s}}\right) \lesssim \eps^2$$
and thus
$$C'(\eps):\approx \left(\norm{H^*}_{\frac{2N}{\alpha}}\norm{K^*}_{\frac{2N}{\alpha}}\right) .$$
In the last term of \eqref{eq_dim_mixed_terms_HK} we use the generalized Young's inequality $ab \leq \frac{\delta}{2} a^2 + \frac{1}{2\delta} b^2$, with
$$\delta :=\eps^2 \left(\norm{H^*}_{\frac{2N}{\alpha}}\norm{K_*}_{\frac{2N}{\alpha+2s}} + \norm{H_*}_{\frac{2N}{\alpha+2s}} \norm{K^*}_{\frac{2N}{\alpha}}\right)^{-1}$$
so that
$$ \left(\norm{H^*}_{\frac{2N}{\alpha}}\norm{K_*}_{\frac{2N}{\alpha+2s}} + \norm{H_*}_{\frac{2N}{\alpha+2s}} \norm{K^*}_{\frac{2N}{\alpha}}\right) \norm{u}_2 \norm{(-\Delta)^{s/2} u}_2 \leq \tfrac{1}{2}\eps^2 \norm{u}_2^2 
+ C''(\eps) \norm{(-\Delta)^{s/2}u}_2^2.$$
Merging the pieces, we have the claim.
\QED

\bigskip

The following technical result can be found in \cite[Lemma 3.5]{GDS}. 

\begin{Lemma}[\cite{GDS}]\label{lem_truncation}
Let $a, b \in \R$, $r \geq 2$ and $k\geq 0$. Set $T_k: \R \to [-k,k]$ the truncation in $k$, that is
$$T_k(t):= \parag{ -k && \quad \hbox{ if $t\leq -k$}, \\ t&& \quad \hbox{ if $t \in (-k,k)$}, \\ k&& \quad \hbox{ if $t \geq k$},}$$
and write $a_k:= T_k(a)$, $b_k:= T_k(b)$. Then
$$\frac{4(r-1)}{r^2} \left(|a_k|^{r/2} - |b_k|^{r/2}\right)^2 \leq (a-b)\left( a_k |a_k|^{r-2} - b_k |b_k|^{r-2}\right).$$
\end{Lemma}

Notice that the (optimal) Sobolev embedding tells us that $H^s(\R^N) \hookrightarrow L^{2^*_s}(\R^N)$. In the following we show that $u$ belongs to some $L^r(\R^N)$ with $r> 2^*_s=\frac{2N}{N-2s}$; notice that we make no use of the Caffarelli-Silvestre $s$-harmonic extension method, and work directly in the fractional framework.

\begin{Proposition}\label{prop_uinLr}
Let $H, K \in L^{\frac{2N}{\alpha}}(\R^N) + L^{\frac{2N}{\alpha+2s}}(\R^N) $. Assume that $u \in H^s(\R^N)$ solves
$$(-\Delta)^s u + u = (I_{\alpha}*(Hu))K, \quad \hbox{ in $\R^N$}$$
in the weak sense. Then
$$u \in L^r(\R^N) \quad \hbox{for all $r \in \left[2, \frac{N}{\alpha} \frac{2N}{N-2s}\right)$}.$$
Moreover, for each of these $r$, we have
$$\norm{u}_r \leq C_r \norm{u}_2$$
with $C_r>0$ not depending on $u$.
\end{Proposition}

\medskip

\claim Proof.
By Lemma \ref{lem_secondMS} there exists $\lambda>0$ (that we can assume large) such that
\begin{equation}\label{eq_H_lambda}
\int_{\R^N} \left(I_{\alpha}*\big(H|u|\big)\right)K|u| \, dx \leq \frac{1}{2} \norm{(-\Delta)^{s/2} u}_2^2 + \frac{\lambda}{2}\norm{u}_2^2.
\end{equation}
Let us set
$$H_n:= H \chi_{\{|H|\leq n\}}, \quad K_n:= K \chi_{\{|K|\leq n\}}, \quad \hbox{ for $n \in \N$}$$
and observe that
$$H_n, \; K_n \in L^{\frac{2N}{\alpha}}(\R^N),$$
$$H_n \to H, \quad K_n \to K \quad \hbox{ almost everywhere, as $n\to +\infty$}$$
and
\begin{equation}\label{eq_H_nH}
|H_n| \leq |H|, \quad |K_n| \leq |K| \quad \hbox{ for every $n \in \N$}.
\end{equation}
We thus define the bilinear form
$$a_n(\varphi, \psi):= \int_{\R^N} (-\Delta)^{s/2} \varphi \, (-\Delta)^{s/2} \psi \, dx + \lambda \int_{\R^N} \varphi \psi \, dx - \int_{\R^N} \left(I_{\alpha}*\big(H_n \varphi\big)\right) K_n \psi \, dx$$
for every $\varphi, \psi \in H^s(\R^N)$. Since, by \eqref{eq_H_nH} and \eqref{eq_H_lambda}, we have
\begin{equation}\label{eq_coercive}
a_n(\varphi, \varphi) \geq \frac{1}{2} \norm{(-\Delta)^{s/2} \varphi}_2^2 + \frac{\lambda}{2} \norm{\varphi}_2^2 \geq \frac{1}{2} \norm{\varphi}_{H^s(\R^N)}^2
\end{equation}
for each $\varphi \in H^s(\R^N)$, we obtain that $a_n$ is coercive. Set
$$f:= (\lambda-1) u \in H^s(\R^N)$$
we obtain by Lax-Milgram theorem that, for each $n \in \N$, there exists a unique $u_n \in H^s(\R^N)$ solution of
$$a_n(u_n, \varphi)= (f, \varphi)_2, \quad \varphi \in H^s(\R^N),$$
that is
\begin{equation}\label{eq_troncat}
(-\Delta)^s u_n + \lambda u_n - \big(I_{\alpha}*(H_n u_n)\big)K_n =(\lambda-1) u, \quad \hbox{ in $\R^N$}
\end{equation}
in the weak sense; moreover the theorem tells us that
$$\norm{u_n}_{H^s} \leq \frac{\norm{f}_2}{1/2}= 2(\lambda-1) \norm{u}_2$$
(since $1/2$ appears as coercivity coefficient in \eqref{eq_coercive}), and thus $u_n$ is bounded. Hence $u_n \wto \bar{u}$ in $H^s(\R^N)$ up to a subsequence for some $\bar{u}$. This means in particular that $u_n \to \bar{u}$ almost everywhere pointwise.
 
 Thus we can pass to the limit in
$$\int_{\R^N} (-\Delta)^{s/2} u_n \, (-\Delta)^{s/2} \varphi \, dx + \lambda \int_{\R^N} u_n \varphi \, dx - \int_{\R^N} \left(I_{\alpha}*\big(H_n u_n\big)\right) K_n \varphi \, dx = (\lambda -1) \int_{\R^N} u \varphi \, dx.$$
We need to check only the Choquard term. 
We first see by the continuous embedding that $u_n \wto \bar{u}$ in $L^q(\R^N)$, for $q \in [2, 2^*_s]$. Split again $H=H^*+ H_*$, $K=K^* + K_*$ and work separately in the four combinations; we assume to work generally with $\tilde{H} \in \{H^*, H_*\}$, $\tilde{H}\in L^{\beta}(\R^N)$ and $\tilde{K} \in \{K^*, K_*\}$, $\tilde{K}\in L^{\gamma}(\R^N)$, where $\beta, \gamma \in \{ \frac{2N}{\alpha}, \frac{2N}{\alpha+2s}\}$. 
Then one can easily prove that $\tilde{H}_n u_n \wto \tilde{H} \bar{u}$ in $L^r(\R^N)$ with $\frac{1}{r} = \frac{1}{\beta} + \frac{1}{q}$.
By the continuity and linearity of the Riesz potential we have $I_{\alpha} * (H_n u_n) \wto I_{\alpha} * (H \bar{u})$ in $L^h(\R^N)$, where $\frac{1}{h}= \frac{1}{r} - \frac{\alpha}{n}$. 
As before, we obtain $\left(I_{\alpha}*\big(H_n u_n\big)\right) K_n \wto \left(I_{\alpha}*\big(H \bar{u}\big)\right) K$ in $L^k(\R^N)$, where $\frac{1}{k} = \frac{1}{\gamma} + \frac{1}{h}$.
Simple computations show that if $\beta=\gamma=\frac{2N}{\alpha}$ and $q=2$, then $k'=2$; if $\beta=\frac{2N}{\alpha}$, $\gamma= \frac{2N}{\alpha+2s}$ (or viceversa) and $q=2$, then $k'=2^*_s$; if $\beta=\gamma=\frac{2N}{\alpha+2s}$ and $q=2^*_s$, then $k'=2^*_s$. 
Therefore $H^s(\R^N) \subset L^{k'}(\R^N)$ and we can pass to the limit in all the four pieces, obtaining
$$\int_{\R^N} \left(I_{\alpha}*\big(H_n u_n\big)\right) K_n \varphi \, dx \to \int_{\R^N} \left(I_{\alpha}*\big(H \bar{u}\big)\right) K \varphi \, dx.$$
Therefore, $\bar{u}$ satisfies
$$(-\Delta)^s \bar{u} + \lambda \bar{u} - \big(I_{\alpha}*(H \bar{u})\big)K =(\lambda-1) u, \quad \hbox{ in $\R^N$}$$
as well as $u$. But we can see this problem, similarly as before, with a Lax-Milgram formulation and obtain the uniqueness of the solution. Thus $\bar{u}=u$ and hence
$$u_n \wto u \quad \hbox{ in $H^s(\R^N)$, as $n\to +\infty$}$$
and almost everywhere pointwise. 
Let now $k\geq 0$ and write
$$u_{n,k}:= T_k(u_n)\in L^2(\R^N) \cap L^{\infty}(\R^N)$$
where $T_k$ is the truncation introduced in Lemma \ref{lem_truncation}. Let $r\geq 2$. We have $|u_{n,k}|^{r/2} \in H^s(\R^N)$, by exploiting \eqref{eq_semin_gagl} and the fact that $h(t):=(T_k(t))^{r/2}$ is a Lipschitz function with $h(0)=0$.
By \eqref{eq_semin_gagl} and by Lemma \ref{lem_truncation} we have
\begin{eqnarray*}
\lefteqn{\frac{4(r-1)}{r^2} \int_{\R^N} |(-\Delta)^{s/2} (|u_{n,k}|^{r/2}) |^2\, dx = C(N,s) \int_{\R^{2N}} \frac{ \frac{4(r-1)}{r^2}\left(|u_{n,k}(x)|^{r/2} - |u_{n,k}(y)|^{r/2}\right)^2}{|x-y|^{N+2s}} \, dx \, dy }\\
&&\leq C(N,s) \int_{\R^{2N}} \frac{\big(u_n(x)-u_n(y)\big)\left(u_{n,k}(x)|u_{n,k}(x)|^{r-2} - u_{n,k}(y) |u_{n,k}(y)|^{r-2}\right)}{|x-y|^{N+2s}} \, dx \, dy .
\end{eqnarray*}
Set
$$\varphi:= u_{n,k}|u_{n,k}|^{r-2}$$
it results that $\varphi \in H^s(\R^N)$, since again $h(t):=T_k(t) |T_k(t)|^{r-2}$ is a Lipschitz function with $h(0)=0$.
Thus we can choose it as a test function in \eqref{eq_troncat} and obtain, by polarizing the identity \eqref{eq_semin_gagl},
\begin{eqnarray*}
\lefteqn{
 \frac{4(r-1)}{r^2} \int_{\R^N} |(-\Delta)^{s/2} (|u_{n,k}|^{r/2}) |^2\, dx \leq C(N,s) \int_{\R^{2N}} \frac{\big(u_n(x)-u_n(y)\big)\left(\varphi(x) - \varphi(y)\right)}{|x-y|^{N+2s}} \, dx \, dy} \\
 &&= - \lambda \int_{\R^N} u_n \varphi \, dx + \int_{\R^N} \left(I_{\alpha}*(H_n u_n)\right) K_n \varphi \, dx + (\lambda -1) \int_{\R^N} u \varphi \, dx
 \end{eqnarray*}
 and since $u_n \varphi \geq |u_{n,k}|^r$ we gain
 \begin{eqnarray}
 \lefteqn{
 \frac{4(r-1)}{r^2} \int_{\R^N} |(-\Delta)^{s/2} (|u_{n,k}|^{r/2}) |^2\, dx \leq } \notag\\
 &&\leq - \lambda \int_{\R^N} |u_{n,k}|^r \, dx+\int_{\R^N} \big(I_{\alpha}*(H_n u_n)\big) K_n \varphi \, dx + (\lambda -1) \int_{\R^N} u \varphi \, dx. \label{eq_dim_4(r-1)}
 \end{eqnarray}
 Focus on the Choquard term on the right-hand side. We have
 \begin{eqnarray}
\lefteqn{\int_{\R^N} \big(I_{\alpha}*(H_n u_n)\big) K_n \varphi \, dx \leq}\\
&\leq& \int_{\R^N} \big(I_{\alpha}*(|H_n| |u_n|\chi_{\{|u_n|\leq k\}})\big) |K_n| |u_{n,k}|^{r-1} \, dx+ \int_{\R^N} \big(I_{\alpha}*(|H_n| |u_n| \chi_{\{|u_n|>k\}})\big) |K_n| |u_{n,k}|^{r-1} \, dx \notag\\
&\leq & \int_{\R^N} \big(I_{\alpha}*(|H_n| |u_{n,k}|)\big) |K_n| |u_{n,k}|^{r-1} \, dx +\int_{\R^N} \big(I_{\alpha}*(|H_n||u_{n}|\chi_{\{|u_n|>k\}})\big) |K_n| |u_n|^{r-1} \, dx \notag\\
&\stackrel{\eqref{eq_H_nH}}\leq & \int_{\R^N} \big(I_{\alpha}*(|H| |u_{n,k}|)\big) |K| |u_{n,k}|^{r-1} \, dx + \int_{\R^N} \big(I_{\alpha}*(|H_n||u_{n}|\chi_{\{|u_n|>k\}})\big) |K_n| |u_n|^{r-1} \, dx \notag\\
& =:& (I)+(II). \label{eq_dim_(I-II)}
 \end{eqnarray}
 
Focus on $(I)$. Consider $r \in [2, \frac{2N}{\alpha})$, so that $\theta:= \frac{2}{r} \in (\frac{\alpha}{N}, 2-\frac{\alpha}{N})$. Choose moreover $v:= |u_{n,k}|^{r/2} \in H^s(\R^N)$ and $\eps^2:= \frac{2(r-1)}{r^2}>0$.
 Thus, observed that if a function belongs to a sum of Lebesgue spaces then its absolute value does the same (\cite[Proposition 2.3]{BPR}), by Lemma \ref{lem_secondMS} we obtain
 \begin{equation}\label{eq_dim_(I)}
 (I) \leq \frac{2(r-1)}{r^2} \norm{(-\Delta)^{s/2}(|u_{n,k}|^{r/2})}_2^2 + C(r) \norm{|u_{n,k}|^{r/2}}_2^2.
 \end{equation}
 Focus on $(II)$. Assuming $r< \min\{\frac{2N}{\alpha}, \frac{2N}{N-2s}\}$, we have $u_n \in L^r(\R^N)$ 
 and $H_n \in L^{\frac{2N}{\alpha}}(\R^N)$, thus
 $$|H_n| |u_n| \in L^{a}(\R^N), \quad \hbox{with $\frac{1}{a} = \frac{\alpha}{2N} + \frac{1}{r}$}$$
 for the H\"older inequality. Similarly
 $$|K_n| |u_n|^{r-1} \in L^{b}(\R^N), \quad \hbox{with $\frac{1}{b} = \frac{\alpha}{2N} + 1-\frac{1}{r}$}.$$
 Thus, since $\frac{1}{a} + \frac{1}{b}= \frac{N+\alpha}{N}$, we have by the Hardy-Littlewood-Sobolev inequality (see Proposition \ref{prop_HLS}) that
$$
\int_{\R^N} \big(I_{\alpha}*(|H_n||u_{n}|\chi_{\{|u_n|>k\}})\big) |K_n| |u_n|^{r-1} \, dx \leq C\left( \int_{\{|u_n|>k\}} \abs{|H_n| |u_n|}^a \, dx\right)^{1/a} \left( \int_{\R^N}\abs{|K_n||u_n|^{r-1}}^{b}\, dx\right)^{1/b} .
$$
 With respect to $k$, the second factor on the right-hand side is bounded, while the first factor goes to zero thanks to the dominated convergence theorem, thus
\begin{equation}\label{eq_dim_(II)}
(II) = o_k(1), \quad \hbox{ as $k\to +\infty$}.
\end{equation}
 Joining \eqref{eq_dim_4(r-1)}, \eqref{eq_dim_(I-II)}, \eqref{eq_dim_(I)}, \eqref{eq_dim_(II)} we obtain
\begin{eqnarray*}
\lefteqn{ \frac{2(r-1)}{r^2} \int_{\R^N} |(-\Delta)^{s/2} (|u_{n,k}|^{r/2}) |^2\, dx \leq} \\
&& \leq - \lambda \int_{\R^N} |u_{n,k}|^r
 \, dx + C(r) \int_{\R^N} |u_{n,k}|^r
 \, dx + (\lambda -1) \int_{\R^N} u \varphi \, dx + o_k(1).
\end{eqnarray*}
 That is, by Sobolev inequality \eqref{eq_embd_homog} 
$$C'(r)\left( \int_{\R^N} |u_{n,k}|^{\frac{r}{2} 2^*_s} \, dx \right)^{2/2^*_s} \leq (C(r)-\lambda) \int_{\R^N} |u_{n,k}|^r \, dx + (\lambda-1)
 \int_{\R^N} |u| \, |u_{n,k}|^{r-1} \, dx + o_k(1).$$
Letting $k\to +\infty$ by the monotone convergence theorem (since $u_{n,k}$ are monotone with respect to $k$ 
and $u_{n,k} \to u_n$ pointwise)
we have
\begin{equation}\label{eq_moser}
C'(r)\left( \int_{\R^N} |u_n|^{\frac{r}{2} 2^*_s} \, dx \right)^{2/2^*_s} \leq (C(r)-\lambda) \int_{\R^N} |u_{n}|^r \, dx + (\lambda-1) \int_{\R^N} |u| \, |u_{n}|^{r-1} \, dx
\end{equation}
and thus $u_n \in L^{\frac{r}{2}2^*_s}(\R^N)$. Notice that $\frac{r}{2} \in \big[1, \min\{\frac{N}{\alpha}, \frac{N}{N-2s}\}\big)$. If $N-2s < \alpha$ we are done. Otherwise, set $r_1:=r$, we can now repeat the argument with
$$r_2 \in \left( \frac{2N}{N-2s}, \min\left\{ \frac{2N}{\alpha}, 2\left(\frac{N}{N-2s}\right)^2\right\}\right).$$
Again, if $\frac{2N}{\alpha} < 2\left(\frac{N}{N-2s}\right)^2$ we are done, otherwise we repeat the argument. Inductively, we have
$$\left(\frac{N}{N-2s}\right)^m \to +\infty, \quad \hbox{as $m\to +\infty$}$$
thus $\frac{2N}{\alpha}<2\left(\frac{N}{N-2s}\right)^m$ after a finite number of steps. For such $r=r_m$, consider again \eqref{eq_moser}: 
by the almost everywhere convergence of $u_n$ to $u$ and Fatou's lemma
\begin{align*}
C''(r) \left(\int_{\R^N} |u|^{\frac{r}{2} 2^*_s}\right)^{2/2^*_s} \, dx & \leq \liminf_{n} C''(r) \left( \int_{\R^N} |u_n|^{\frac{r}{2} 2^*_s} \, dx \right)^{2/2^*_s} \\ 
&\leq \liminf_n \left( (C(r)-\lambda) \int_{\R^N} |u_{n}|^r \, dx + (\lambda-1) \int_{\R^N} |u| \, |u_{n}|^{r-1} \, dx\right) \\
&\leq (C(r)-\lambda ) \limsup_n\int_{\R^N} |u_{n}|^r \, dx + (\lambda-1) \limsup_n \int_{\R^N} |u| \, |u_{n}|^{r-1} \, dx.
\end{align*}
Being $u_n$ equibounded in $H^s(\R^N)$ and thus in $L^{2^*_s}(\R^N)$, by the iteration argument we have that it is equibounded also in $L^r(\R^N)$; in particular, the bound is given by $\norm{u}_2$ times a constant $C(r)$. Thus the right-hand side is a finite quantity, and we gain $u \in L^{\frac{r}{2}2^*_s}(\R^N)$, which is the claim.
\QED

\bigskip

The following Lemma states that $I_{\alpha}*g \in L^{\infty}(\R^N)$ whenever $g$ lies in $L^q(\R^N)$ with $q$ in a neighborhood of $\frac{N}{\alpha}$ (in particular, it generalizes Proposition \ref{prop_HLS} to the case $h=\infty$ and $r \approx \frac{N}{\alpha}$).

In addition, it shows the decay at infinity of the Riesz potential, which will be useful in Section \ref{sec_asymptotic}.

\begin{Proposition}\label{prop_conv_C0}
Assume that \textnormal{(f1)-(f2)} hold. Let $u\in H^s(\R^N)$ be a solution of \eqref{problem_x}. Then 
$u\in L^q(\R^N)$ for $q \in \big[2,\frac{N}{\alpha} \frac{2N}{N-2s}\big)$, and
$$I_{\alpha} * F(u) \in C_0(\R^N),$$
that is, continuous and zero at infinity. 
In particular, 
$$I_{\alpha} * F(u) \in L^{\infty}(\R^N)$$
and
$$\big(I_{\alpha} * F(u)\big)(x) \to 0 \quad \hbox{as $|x| \to +\infty$}.$$
\end{Proposition}

\claim Proof.
We first check to be in the assumptions of Proposition \ref{prop_uinLr}. Indeed, by (f1)-(f2) and the fact that $u\in H^s(\R^N)\subset L^2(\R^N) \cap L^{2^*_s}(\R^N)$ we obtain that
$$H:= \frac{F(u)}{u}, \quad K:=f(u)$$
lie in $L^{\frac{2N}{\alpha}}(\R^N) + L^{\frac{2N}{\alpha+2s}}(\R^N)$, since bounded by functions in this sum space (see e.g. \cite[Proposition 2.3]{BPR}). 
Now by Proposition \ref{prop_uinLr} we have $u \in L^q(\R^N)$ for $q \in [2, \frac{N}{\alpha} \frac{2N}{N-2s})$. 

To gain the information on the convolution, we want to use Young's Theorem, which states that if $g, h$ belong to two Lebesgue spaces with conjugate (finite) indexes, then $g*h \in C_0(\R^N)$. 
We first split
$$I_{\alpha}*F(u) = (I_{\alpha}\chi_{B_1})*F(u) + (I_{\alpha}\chi_{B_1^c})*F(u)$$
where
$$ I_{\alpha}\chi_{B_1} \in L^{r_1}(\R^N), \quad \hbox{ for $r_1 \in [1, \frac{N}{N-\alpha})$},$$
$$ I_{\alpha}\chi_{B_1^c} \in L^{r_2}(\R^N), \quad \hbox{ for $r_2 \in (\frac{N}{N-\alpha}, \infty]$}.$$
We need to show that $F(u) \in L^{q_1}(\R^N)\cap L^{q_2}(\R^N)$ for some $q_i$ satisfying
$$\frac{1}{q_i} + \frac{1}{r_i} = 1, \quad i=1,2$$
that is
$$\frac{q_1}{q_1-1} \in \left[1, \frac{N}{N-\alpha}\right), \quad \frac{q_2}{q_2-1}\in \left(\frac{N}{N-\alpha}, \infty\right]$$
or equivalently $q_2 < \frac{N}{\alpha} < q_1$.
Recall that
$$|F(u)| \leq C\left(|u|^{\frac{N+\alpha}{N}} + |u|^{\frac{N+\alpha}{N-2s}}\right).$$
Note that $u \in L^q(\R^N)$ for $q \in [2, \frac{N}{\alpha} \frac{2N}{N-2s})$ implies
$$|u|^{\frac{N+\alpha}{N}} , |u|^{\frac{N+\alpha}{N-2s}} \in L^{q_1}(\R^N) \cap L^{q_2}(\R^N)$$
for some $q_2<\frac{N}{\alpha}<q_1$. Thus we have the claim.
\QED

\bigskip

Once obtained the boundedness of the Choquard term, we can finally gain the boundedness of the solution.

\begin{Proposition}\label{prop_u_bounded}
Assume that \textnormal{(f1)-(f2)} hold. Let $u\in H^s(\R^N)$ be a positive solution of \eqref{problem_x}. Then $u\in L^{\infty}(\R^N)$.
\end{Proposition}

\claim Proof.
By Lemma \ref{prop_conv_C0} 
we obtain
$$a:= I_{\alpha}*F(u) \in L^{\infty}(\R^N).$$
Thus $u$ satisfies the following nonautonomous problem, with a local nonlinearity
$$(-\Delta)^{s/2} u + \mu u = a(x) f(u), \quad \hbox{ in $\R^N$}$$
with $a$ bounded. In particular
$$(-\Delta)^{s/2} u = g(x,u):=- \mu u + a(x) f(u), \quad \hbox{ in $\R^N$}$$
where
$$|g(x,t)| \leq \mu |t| + C \norm{a}_{\infty} \left(|t|^{\frac{\alpha}{N}} + |t|^{\frac{\alpha+2s}{N-2s}}\right).$$
Set $\gamma:= \max\{1, \frac{\alpha+2s}{N-2s}\} \in [1, 2^*_s)$, we thus have
$$|g(x,t)| \leq C(1 + |t|^{\gamma}).$$
Hence we are in the assumptions of \cite[Proposition 5.1.1]{DMV} and we can conclude. 
\QED

\bigskip

We observe that a direct proof of the boundedness for generally signed solutions, but assuming also (f5), can be found in Appendix \ref{sec_bound_sign}. 

\smallskip

Gained the boundedness of the solutions, we obtain also some additional regularity, which will be implemented in some bootstrap argument for the $L^1$-summability.

\begin{Proposition}\label{prop_u_holder}
Assume that \textnormal{(f1)-(f2)} hold. Let $u\in H^s(\R^N)\cap L^{\infty}(\R^N)$ be a weak solution of \eqref{problem_x}. Then $u \in H^{2s}(\R^N) \cap C^{0,\gamma}(\R^N)$ for any $\gamma \in (0, \min\{1,2s\})$. Moreover $u$ satisfies \eqref{problem_x} almost everywhere.
\end{Proposition}

\claim Proof.
By Proposition \ref{prop_u_bounded}, Proposition \ref{prop_conv_C0} and (f2) we have that $u\in L^{\infty}(\R^N)$ satisfies
$$(-\Delta)^s u = g \in L^{\infty}(\R^N)$$
where $g(x):=- \mu u(x) + (I_{\alpha}*F(u))(x) f(u(x))$. 
We prove first that $u\in H^{2s}(\R^N)$. Indeed, we already know that $f(u)$, $F(u)$ and $I_{\alpha}*F(u)$ belong to $L^{\infty}(\R^N)$. By Remark \ref{rem_buona_posit}, we obtain
$$f(u) \in L^{\frac{2N}{\alpha+2s}}(\R^N)\cap L^{\infty}(\R^N), \quad F(u) \in L^{\frac{2N}{N+\alpha}}(\R^N) \cap L^{\infty}(\R^N),$$
$$I_{\alpha}*F(u) \in L^{\frac{2N}{N-2s}}(\R^N) \cap L^{\infty}(\R^N), \quad (I_{\alpha}*F(u))f(u) \in L^2(\R^N)\cap L^{\infty}(\R^N).$$
In particular,
$$g:=(I_{\alpha}*F(u))f(u)-\mu u \in L^2(\R^N).$$
Since $u$ is a weak solution, we have, fixed $\varphi \in H^s(\R^N)$,
\begin{equation}\label{eq_dim_weak_form}
\int_{\R^N} (-\Delta)^{s/2} u \, (-\Delta)^{s/2} \varphi \, dx = \int_{\R^N} g \, \varphi \, dx.
\end{equation}
Since $g \in L^2(\R^N)$, we can apply Plancharel theorem and obtain
\begin{equation}\label{eq_dim_plancharel}
\int_{\R^N} |\xi|^{2s} \widehat{u} \, \widehat{\varphi} \, d\xi = \int_{\R^N} \widehat{g} \, \widehat{\varphi} \, d \xi.
\end{equation}
Since $H^s(\R^N) = \mc{F}(H^s(\R^N))$ and $\varphi$ is arbitrary, we gain
$$|\xi|^{2s} \widehat{u} = \widehat{g} \in L^2(\R^N).$$
By definition, we obtain $u \in H^{2s}(\R^N)$, which concludes the proof.
Observe moreover that $\mc{F}^{-1}\big((1+|\xi|^{2s})\widehat{u}\big) = u +g \in L^2(\R^N) \cap L^{\infty}(\R^N)$, thus by definition $u \in H^{2s}(\R^N) \cap W^{2s, \infty}(\R^N)$. By the embedding \cite[Theorem 3.2]{FQT} we obtain $u \in C^{0,\gamma}(\R^N)$ if $2s<1$ and $\gamma \in (0, 2s)$, while $u \in C^{1,\gamma}(\R^N)$ if $2s>1$ and $\gamma \in (0,2s-1)$ 
(see also \cite[Proposition 2.9]{Sil0}).

It remains to show that $u$ is an almost everywhere pointwise solution. Thanks to the fact that $u\in H^{2s}(\R^N)$, we use again \eqref{eq_dim_plancharel}, where we can apply Plancharel theorem (that is, we are integrating by parts \eqref{eq_dim_weak_form}) and thus
$$\int_{\R^N}(-\Delta)^s u \, \varphi \, dx = \int_{\R^N} g \, \varphi \, dx.$$
Since $\varphi \in H^s(\R^N)$ is arbitrary, we obtain
$$(-\Delta)^s u = g \quad \hbox{ almost everywhere.}$$
This concludes the proof.
\QED

\bigskip

We observe, by the proof, that if $s\in (\tfrac{1}{2}, 1)$, then $u \in C^{1,\gamma}(\R^N)$ for any $\gamma \in (0, 2s-1)$, and $u$ is a classical solution, with $(-\Delta)^s u \in C(\R^N)$ and equation \eqref{problem_x} satisfied pointwise.

\smallskip

We end this section by dealing with the summability of $u$ in Lebesgue spaces $L^r(\R^N)$ for $r<2$.

\begin{Remark}\label{rem_L1}
We start noticing that, if a solution $u$ belongs to some $L^q(\R^N)$ with $q<2$, then $u\in L^1(\R^N)$. 
Assume thus $q \in (1,2)$ and let $u \in L^q(\R^N) \cap L^{\infty}(\R^N)$, then we have
$$f(u) \in L^{\frac{qN}{\alpha}}(\R^N)\cap L^{\infty}(\R^N), \quad F(u) \in L^{\frac{qN}{N+\alpha}}(\R^N) \cap L^{\infty}(\R^N),$$
$$I_{\alpha}* F(u) \in L^{\frac{qN}{N+ \alpha(1-q)}}(\R^N) \cap L^{\infty}(\R^N), \quad (I_{\alpha}*F(u)) f(u) \in L^{\frac{qN}{N + \alpha(2-q)}}(\R^N) \cap L^{\infty}(\R^N).$$
Thanks to Proposition \ref{prop_u_holder}, $u$ satisfies \eqref{problem_x} almost everywhere, thus we have
$$\mc{F}^{-1}\big((|\xi|^{2s} + \mu)\, \widehat{u}\big) = (-\Delta)^s u + \mu u = (I_{\alpha}*F(u))f(u) \in L^{\frac{qN}{N + \alpha(2-q)}}(\R^N)$$
which equivalently means that the Bessel operator verifies
$$\mc{F}^{-1}\big((|\xi|^{2} + 1)^s \,\widehat{u}\big) \in L^{\frac{qN}{N + \alpha(2-q)}}(\R^N).$$
Thus by \cite[Theorem 1.2.4]{AH} we obtain that $u$ itself lies in the same Lebesgue space, that is
$$u \in L^{\frac{qN}{N + \alpha(2-q)}}(\R^N).$$
If $\frac{qN}{N + \alpha(2-q)}<1$, we mean that $(I_{\alpha}*F(u)) f(u) \in L^1(\R^N) \cap L^{\infty}(\R^N)$, and thus $u \in L^1(\R^N) \cap L^{\infty}(\R^N)$. We convey this when we deal with exponents less than $1$.

If $q <2$, then
$$\frac{qN}{N + \alpha(2-q)}<q$$
and we can implement a bootstrap argument to gain $u \in L^1(\R^N)$. More precisely
$$\parag{ &q_0\in [1, 2)& \\ &q_{n+1} = \frac{q_n N}{N + \alpha(2-q_n)}& }$$
where $q_n \to 0$ (but we stop at $1$). 
Thus, in order to implement the argument, we need to show that $u\in L^q(\R^N)$ for some $q<2$.
\end{Remark}

We show now that $u\in L^1(\R^N)$. 
It is easy to see that, if the problem is (strictly) not lower-critical, i.e. (f2) holds together with
	$$	\lim_{t\to 0}{F(t)\over \abs t^\beta}=0
	$$
for some $\beta\in ({N+\alpha\over N}, {N+\alpha\over N-2s})$, then $u\in L^1(\R^N)$. 
Indeed $u\in H^s(\R^N)\cap L^\infty(\R^N)\subset
L^2(\R^N)\cap L^\infty(\R^N)$ and
	$$	(I_\alpha*F(u))f(u) \in L^q(\R^N),
	$$
where ${1\over q}={\beta\over 2}-{\alpha\over 2N}$; noticed that $q<2$, we can implement the bootstrap argument of Remark \ref{rem_L1}.

We will show that the same conclusion can be reached by assuming only (f2).

\begin{Proposition}\label{prop_u_L1}
Assume that \textnormal{(f1)}-\textnormal{(f2)} hold. Let $u\in H^s(\R^N)\cap L^{\infty}(\R^N)$ be a weak solution of \eqref{problem_x}. 
Then $u \in L^1(\R^N)$.
\end{Proposition}

\medskip

\claim Proof of Proposition \ref{prop_u_L1}.
For a given solution $u\in H^s(\R^N)\cap L^\infty(\R^N)$ we
set again
 $$ H:={F(u)\over u}, \quad K:=f(u).
 $$
Since $u\in L^2(\R^N)\cap L^\infty(\R^N)$, by (f2) we have
$H$, $K\in L^{2N\over \alpha}(\R^N)$. For $n\in\N$, we set
 $$ H_n:=H\chi_{\{ |x|\geq n\} }.
 $$
Then we have
 \begin{equation}\label{kt.a}
 \norm{H_n}_{2N\over \alpha}\to 0 \quad
 \hbox{as}\ n\to\infty.
 \end{equation}
Since $\hbox{supp} (H-H_n)\subset \big\{ |x|\leq n\big\}$ is a
bounded set, we have for any $\beta \in [1,{2N\over \alpha}]$ 
 \begin{equation}\label{kt.b}
 H-H_n \in L^\beta(\R^N) \quad \hbox{for all}\ n\in\N.
 \end{equation}
We write our equation \eqref{problem_x} as
 $$ (-\Delta)^s u+\mu u =(I_\alpha*H_nu)K +R_n \quad \hbox{in $\mathbb{R}^N$},
 $$
where we introduced the function $R_n$ by
 $$ R_n:= (I_\alpha*(H-H_n)u)K.
 $$
Now we consider the following linear equation:
 \begin{equation}\label{kt.c}
 (-\Delta)^s v+\mu v =(I_\alpha*H_nv)K +R_n \quad \hbox{in $\mathbb{R}^N$}.
 \end{equation}
We have the following facts:
\begin{itemize}
\item[(i)] The given solution $u$ solves \eqref{kt.c}.
\item[(ii)] By the property \eqref{kt.b} with $\beta \in (\frac{2N}{N+\alpha},{2N\over \alpha})$, there exists $q_1
\in (1,2)$, namely ${1\over q_1}={1\over\beta}+\half-{\alpha\over 2N}$, 
such that $R_n\in L^{q_1}(\R^N)\cap L^2(\R^N)$.
\item[(iii)] By the property \eqref{kt.a}, for any $r\in (\frac{2N}{2N-\alpha}, 2] \subset (1,2]$
 $$ v\in L^r(\R^N) \mapsto A_n(v):=(I_\alpha*H_n v)K\in L^r(\R^N)
 $$
is well-defined and verifies
 \begin{equation}\label{kt.d}
 \norm{A_n(v)}_r \leq C_{r,n}\norm v_r.
 \end{equation}
Here $C_{r,n}$ satisfies $C_{r,n}\to 0$ as $n\to\infty$.
\end{itemize}
We show only (iii). Since $v\in L^r(\R^N)$, by Hardy-Littlewood-Sobolev inequality and H\"older inequality we obtain
 $$ \norm{A_n(v)}_r \leq C_r\norm{H_n}_{2N\over \alpha}
 \norm K_{2N\over \alpha} \norm v_r,
 $$
where $C_r>0$ is independent of $n$, $v$. Thus by \eqref{kt.a} we have
$C_{r,n}:=C_r\norm{H_n}_{2N\over \alpha}\norm K_{2N\over \alpha}\to 0$ 
as $n\to\infty$.

Now we show $u\in L^{q_1}(\R^N)$, where $q_1\in (1,2)$ is
given in (ii).
Since $((-\Delta)^s+\mu)^{-1}:\, L^r(\R^N)\to L^r(\R^N)$ is
a bounded linear operator for $r\in(1,2]$ (see \cite[Theorem 
1.2.4]{AH}), \eqref{kt.c} can be rewritten as
 $$ v=T_n(v),
 $$
where
 $$ T_n(v):=((-\Delta)^s+\mu)^{-1}\big(A_n(v)+R_n\big).
 $$
By choosing $\beta \in (2, \frac{2N}{\alpha})$ we have $q_1 \in (\frac{2N}{2N-\alpha},2)\subset(1,2)$, thus we observe that for $n$ large, $T_n$ is a contraction in $L^2(\R^N)$ and $L^{q_1}(\R^N)$. We fix such an $n$.

Since $T_n$ is a contraction in $L^2(\R^N)$, we can see that
$u\in H^s(\R^N)$ is a unique fixed point of $T_n$.
In particular, we have
 $$ u = \lim_{k\to \infty} T_n^k (0) \quad \hbox{in} \ L^2(\R^N).
 $$
On the the other hand, since $T_n$ is a contraction in $L^{q_1}(\R^N)$,
$(T_n^k (0))_{k=1}^\infty$ also converges in $L^{q_1}(\R^N)$. Thus the
limit $u$ belongs to $L^{q_1}(\R^N)$.

Since $q_1<2$ we can use the bootstrap argument of Remark \ref{rem_L1} to get $u\in L^1(\R^N)$, and reach the claim.
\QED

\section{Asymptotic decay}\label{sec_asymptotic}

We prove now the polynomial decay of the solutions. We start from two standard lemmas, whose proofs can be found for instance in \cite[Lemma A.1 and Lemma A.3]{CG0}.

\begin{Lemma}[Maximum Principle]\label{lem_comp_prin}
Let $\Sigma \subset \R^N$, possibly unbounded, and let $u\in H^s(\R^N)$ be a weak subsolution of
$$(-\Delta)^s u + a u \leq 0 \quad \hbox{in $\mathbb{R}^N\setminus \Sigma$}$$
with $a>0$, in the sense that
$$\int_{\R^N} (-\Delta)^{s/2} u \,(-\Delta)^{s/2} \varphi \, dx+ a \int_{\R^N} u \varphi \, dx\leq 0$$
for every positive $\varphi \in H^s(\R^N)$ with $\supp(\varphi) \subset \R^N \setminus \Sigma$. 
Assume moreover that
$$u\leq 0, \quad \textit{ for a.e. $x \in \Sigma$}.$$
Then
\begin{equation}\label{eq_dis_comp_prin}
u\leq 0, \quad \textit{ for a.e. $x \in \R^N$}.
\end{equation}
\end{Lemma}

\begin{Lemma}[Comparison function]\label{lem_esist_sol_part}
Let $b>0$. Then there exists a strictly positive continuous function $W\in H^s(\R^N)$ such that, for some positive constants $C', C''$ (depending on $b$), it verifies 
$$(-\Delta)^s W + b W = 0 \quad \hbox{in $\mathbb{R}^N\setminus B_{r}$}$$
pointwise, with $r:= b^{-1/2s}$, and
\begin{equation}\label{eq_stima_fun_confr}
\frac{C'}{|x|^{N+2s}}<W(x)< \frac{C''}{|x|^{N+2s}}, \quad \textit{ for $|x|>2 r$}.
\end{equation}
\end{Lemma}

We show first some conditions which imply the decay at infinity of the solutions.

\begin{Lemma}\label{lem_u_to0}
Assume that \textnormal{(f1)-(f2)} hold. Let $u$ be a weak solution of \eqref{problem_x}. Assume 
$$u \in L^{\frac{N}{2s}}(\R^N)\cap L^{\infty}(\R^N)$$
and
$$(I_{\alpha}*F(u))f(u)\in L^{\frac{N}{2s}}(\R^N)\cap L^{\infty}(\R^N).$$
 Then we have
\begin{equation}\label{eq_conv_unif_0_veps2}
u(x) \to 0 \quad \hbox{as $|x|\to +\infty$}.
\end{equation}
\end{Lemma}

\claim Proof.
Being $u$ solution of
$$(-\Delta)^s u + u = (1-\mu) u + \big(I_{\alpha}*F(u)\big) f(u) =: \chi \quad \hbox{in $\mathbb{R}^N$},$$
where $\chi \in L^{\frac{N}{2s}}(\R^N) \cap L^{\infty}(\R^N)$, 
we have the representation formula
$$ u= \mc{K} * \chi $$
where $\mc{K}$ is the Bessel kernel; we recall that $\mc{K}$ is positive, it satisfies $\mc{K}(x) \leq \frac{C}{|x|^{N+2s}}$ for $|x| \geq 1$ and $\mc{K} \in L^q(\R^N)$ for $q \in [1, 1 + \tfrac{2s}{N-2s})$ (see \cite[page 1241 and Theorem 3.3]{FQT}). Let us fix $\eta>0$; 
we have, for $x \in \R^N$,
\begin{align*}
u(x) =& \int_{\R^N} \mc{K}(x-y) \chi(y) dy \\
=& \int_{|x-y|\geq 1/\eta} \mc{K}(x-y) \chi(y)dy +\int_{|x-y|< 1/\eta} \mc{K}(x-y) \chi(y)dy.
\end{align*}
As regards the first piece
$$ \int_{|x-y|\geq 1/\eta} \mc{K}(x-y) \chi(y)dy \leq \norm{\chi}_{\infty} \int_{|x-y|\geq 1/\eta} \frac{C}{|x-y|^{N+2s}} dy \leq C \eta^{2s}$$
while for the second piece, fixed a whatever $q \in (1, 1 + \tfrac{2s}{N-2s})$ and its conjugate exponent $q'> \frac{N}{2s}$, we have by H\"older inequality
$$
\int_{|x-y|< 1/\eta} \mc{K}(x-y) \chi(y)dy \leq \norm{\mc{K}}_q \norm{\chi}_{L^{q'}(B_{1/\eta}(x))}\\
$$
where the second factor can be made small for $|x| \gg 0$. 
Joining the pieces, we have \eqref{eq_conv_unif_0_veps2}. 
\QED

\bigskip

We observe that the assumptions of the Lemma are fulfilled by assuming that $u$ is bounded thanks to Proposition \ref{prop_u_L1}. 
We are now ready to prove the polynomial decay of the solutions.

\medskip

\claim Conclusion of the proof of Theorem \ref{th_INT_decay}.
Observe that, by (f5) and Lemma \ref{lem_u_to0}, we have 
\begin{equation}\label{eq_dim_V_lim}
\frac{f(u)}{u} \in L^{\infty}(\R^N).
\end{equation}
Thus we obtain, by applying Proposition \ref{prop_conv_C0}, that
\begin{equation}\label{eq_conv_inf_block}
(I_{\alpha}*F(u))(x) \frac{f(u(x))}{u(x)} \to 0 \quad \hbox{ as $|x| \to + \infty$}. 
\end{equation}
Thus, by \eqref{eq_conv_inf_block} and the positivity of $u$, we have for some $R'\gg 0$
$$(-\Delta)^s u + \tfrac{1}{2} \mu u = (I_{\alpha}*F(u))f(u) - \tfrac{1}{2} \mu u = \left( (I_{\alpha}*F(u))\tfrac{f(u)}{u} - \tfrac{1}{2} \mu \right) u \leq 0 \quad \hbox{in $\mathbb{R}^N\setminus B_{R'}$}.$$
Similarly
$$(-\Delta)^s u + \tfrac{3}{2} \mu u = (I_{\alpha}*F(u))f(u) + \tfrac{1}{2} \mu u = \left( (I_{\alpha}*F(u))\tfrac{f(u)}{u} + \tfrac{1}{2} \mu \right) u \geq 0 \quad \hbox{in $\mathbb{R}^N\setminus B_{R'}$}.$$
Notice that we always intend differential inequalities in the weak sense, that is tested with functions in $H^s(\R^N)$ with supports contained in the reference domain (e.g. $\R^N \setminus B_{R'}$). 

In addition, by Lemma \ref{lem_esist_sol_part} we have that there exist two positive functions $\underline{W}'$, $\overline{W}'$ and three positive constants $R''$, $C'$ and $C''$ depending only on $\mu$, such that
$$ \parag{
& (-\Delta)^s \underline{W}' + \frac{3}{2}\mu \, \underline{W}' = 0 \quad \hbox{in $\mathbb{R}^N \setminus B_{R''}$},& \\ 
&\frac{C'}{|x|^{N+2s}}< \underline{W}' (x), \quad \textnormal{ for $|x|>2R''$}.&}$$
and
$$ \parag{& (-\Delta)^s \overline{W}' + \frac{1}{2}\mu \, \overline{W}' = 0 \quad \hbox{in $\mathbb{R}^N \setminus B_{R''}$},& \\ 
& \overline{W}'(x) < \frac{C''}{|x|^{N+2s}}, \quad \textnormal{ for $|x|>2R''$}.&}$$
Set $R:=\max\{ R', 2R''\}$. Let $\underline{C}_1$ and $\overline{C}_1$ be some lower and upper bounds for $u$ on $B_R$, $\underline{C}_2:=\min_{B_R} \overline{W}'$ and $\overline{C}_2:= \max_{B_R} \underline{W}'$, all strictly positive. Define
$$\underline{W}:= \underline{C}_1 \overline{C}_2 ^{-1} \underline{W}', \quad \overline{W}:= \overline{C}_1 \underline{C}_2^{-1} \overline{W}'$$
so that
$$\underline{W}(x)\leq u(x) \leq \overline{W}(x), \quad \textnormal{ for $|x|\leq R$}.$$
Thanks to the comparison principle in Lemma \ref{lem_comp_prin}, and redefining $C'$ and $C''$, we obtain
$$ \frac{C'}{|x|^{N+2s}} <\underline{W}(x) \leq u(x) \leq \overline{W}(x) < \frac{C''}{|x|^{N+2s}}, \quad \textnormal{ for $|x|>R$}.$$
By the boundedness of $u$, we obtain the claim.
\QED

\bigskip

We see that, for non sublinear $f$ (that is, (f5)), the decay is essentially given by the fractional operator. It is important to remark that, contrary to the limiting local case $s=1$ (see \cite{MS0}), the Choquard term in case of linear $f$ does not affect the decay of the solution. 

\begin{Remark}
We observe that the conclusion of the proof of Theorem \ref{th_INT_decay} can be substituted by exploiting a result in \cite{FLS}. Indeed write $V:= -(I_{\alpha}*F(u)) \frac{f(u)}{u}$, which is bounded and zero at infinity as observed in \eqref{eq_dim_V_lim}-\eqref{eq_conv_inf_block}, and gain
$$(-\Delta)^s u + V(x) u = - \mu u \quad \hbox{in $\mathbb{R}^N$}. $$
Up to dividing for $\norm{u}_2$, we may assume $\norm{u}_2=1$. Thus we are in the assumptions of \cite[Lemma C.2]{FLS} and obtain, even for changing-sign solutions of \eqref{problem_x},
$$|u(x)| \leq \frac{C_1}{(1 + |x|^2)^{\frac{N+2s}{2}}}$$
together with
$$|u(x)| = \frac{C_2}{|x|^{N+2s}} + o\left( \frac{1}{|x|^{N+2s}}\right) \quad \hbox{ as $|x| \to +\infty$}$$
for some $C_1, C_2>0$.
\end{Remark}

\appendix

\section{Boundedness of signed solutions} \label{sec_bound_sign}

In order to achieve the boundedness of general signed solution, we ask in addition that $f$ satisfies (f5). We adapt some argument from \cite[Proposition 2.3]{Gal0}, giving here the details for the reader's convenience. 

\begin{Proposition}\label{eq_prop_u_Linf}
Let $u\in H^s(\R^N)$ be a weak subsolution of
$$(-\Delta)^s u \leq g(x,u) \quad \hbox{in $\mathbb{R}^N$}$$
with
$$|g(x,t)| \leq C\big(|t| 
+ |t|^{2^*_s-1}\big).$$
Then $u \in L^{\infty}(\R^N)$.
\end{Proposition}

In particular this apply to \eqref{problem_x} with
$$g(x,u) := (I_{\alpha}*F(u))f(u) - \mu u,$$
whenever (f5) holds (together with (f1)-(f2)), thanks to Proposition \ref{prop_conv_C0}.
\medskip

\claim Proof.
We already know that $ u \in L^{2^*_s}(\R^N)$. Let us introduce $\gamma >1$, to be fixed, and an arbitrary $T>0$, and set a \emph{$\gamma$-linear (positive) truncation at $T$}
$$h(t)\equiv h_{T, \gamma}(t):= \parag{ &0& \quad \hbox{ if $t \leq 0$}, \\ &t^{\gamma}& \quad \hbox{if $t \in (0, T]$}, \\ &\gamma T^{\gamma-1} t - (\gamma -1) T^{\gamma}& \quad \hbox{ if $t > T$}.}$$
We have that $h\in C^1(\R)\cap W^{1,\infty}(\R)$, it is positive (increasing and convex), zero on the negative halfline, and by direct computations it satisfies the following properties
\begin{equation}\label{eq_dim_h_1}
0 \leq h(t)\leq |t|^{\gamma}, \quad t \in \R,
\end{equation}
\begin{equation}\label{eq_dim_h_2}
0\leq t h'(t) \leq \gamma h(t), \quad t \in \R,
\end{equation}
\begin{equation}\label{eq_dim_h_3}
\lim_{T\to +\infty} h_{T, \gamma}(t) = t^{\gamma}, \quad t \geq 0.
\end{equation}
The goal is to estimate $\norm{h( u)}_{2^*_s}$ and give thus a bound of $ u$ in $L^{2^*_s \gamma}(\R^N)$, where $2^*_s \gamma > 2^*_s$. In order to handle the weak formulation of the notion of solution we introduce 
$$\tilde{h}(t):= \int_0^t (h'(r))^2 \, dr, \quad t \in \R$$
and observe that $\tilde{h}\in C^1(\R) \cap W^{1, \infty}(\R)$ is positive, increasing, convex and zero on the negative halfline. In particular
\begin{equation}\label{eq_dim_htild_0}
\tilde{h}'(t) = (h'(t))^2, \quad t \in \R
\end{equation}
by definition and
\begin{equation}\label{eq_dim_htild_1}
\tilde{h}(t)-\tilde{h}(r) \leq \tilde{h}'(t) (t-r), \quad t, \,r \in \R
\end{equation}
by convexity, and we gain also the Lipschitz continuity
\begin{equation*}\label{eq_dim_htild_abs_1}
|\tilde{h}(t)-\tilde{h}(r)| \leq \norm{\tilde{h}'}_{\infty} |t-r|, \quad t, \,r \in \R.
\end{equation*}
 Combining the definition of $\tilde{h}$, \eqref{eq_dim_h_2} and \eqref{eq_dim_h_1} we obtain
\begin{equation}\label{eq_dim_htild_2}
0 \leq \tilde{h}(t) \leq \norm{h'}_{\infty} |t|^{\gamma}, \quad t \in \R.
\end{equation}
Finally, by a direct application of Jensen inequality we gain
\begin{equation}\label{eq_dim_h_htild}
|h(t)-h(r)|^2 \leq \big( \tilde{h}(t) - \tilde{h}(r)\big) (t-r) , \quad t,\, r \in \R.
\end{equation}
We observe that $\tilde{h}( u) \in H^s(\R^N)$ since $\tilde{h}$ is Lipschitz continuous and $\tilde{h}(0)=0$;
moreover, since $2^*_s$ is the best summability exponent, if we assume 
\begin{equation}\label{eq_dim_stima_gamma}
1<\gamma \leq \frac{2^*_s}{2}
\end{equation}
by \eqref{eq_dim_htild_2} we obtain also 
$$\tilde{h}( u) \leq \norm{h'}_{\infty} |u|^{\gamma} \in L^2(\R^N).$$

We use now the embedding \eqref{eq_embd_homog} and combine \eqref{eq_semin_gagl}, \eqref{eq_dim_h_htild} and the polarized version of \eqref{eq_semin_gagl} to obtain
\begin{align*}
\norm{h( u)}_{2^*_s}^2 &\leq \mc{S}^{-1} \norm{(-\Delta)^{s/2} h( u)}_2^2 \\
&= (C'(N,s))^{-1} \mc{S}^{-1} \int_{\R^{2N}} \frac{|h( u(x))-h( u(y))|^2}{|x-y|^{N+2s}} \, dx \, dy \\
&\leq (C'(N,s))^{-1} \mc{S}^{-1} \int_{\R^{2N}} \frac{\big(\tilde{h}( u(x))- \tilde{h}( u(y))\big)\big( u(x)- u(y)\big)}{|x-y|^{N+2s}} \, dx \, dy \\
&= \mc{S}^{-1} \int_{\R^N} (-\Delta)^{s/2} u \, (-\Delta)^{s/2} \tilde{h}( u) \, dx.
\end{align*}
Since $\tilde{h}( u) \in H^s(\R^N)$ we can choose it as a test function in the equation and gain
$$\norm{h( u)}_{2^*_s}^2 \leq \mc{S}^{-1} \int_{\R^N} g(x,u) 
 \tilde{h}( u) \, dx.$$
By the assumptions on $g$ 
and the positivity of $\tilde{h}( u)$ we obtain
\begin{equation*}\label{eq_posit_to_negat}
\norm{h(u)}_{2^*_s}^2 \leq \mc{S}^{-1} \int_{\R^N} |g(x,u)| \tilde{h}( u) \, dx \leq C \mc{S}^{-1} \int_{\R^N}\big(|u|+ |u|^{2^*_s-1} \big)\tilde{h}( u) \, dx.
\end{equation*}
Since $h(u)$ and $\tilde{h}(u)$ are zero when $u$ is negative, we obtain
$$\norm{h(u_+)}_{2^*_s}^2 \leq C \mc{S}^{-1} \int_{\R^N}\big(u_++ u_+^{2^*_s-1} \big)\tilde{h}( u_+) \, dx.$$
Now we use \eqref{eq_dim_htild_1} (with $r=0$), \eqref{eq_dim_htild_0}, and \eqref{eq_dim_h_2}
\begin{eqnarray}
\lefteqn{ \norm{h(u_+)}_{2^*_s}^2 \leq C \mc{S}^{-1} \int_{\R^N} \big(u_++ u_+^{2^*_s-1} \big) u_+ \tilde{h}'( u_+) \, dx} \nonumber \\
&\leq& C \mc{S}^{-1} \int_{\R^N} \big(u_++ u_+^{2^*_s-1} \big) u_+ (h'(u_+))^2 \, dx \leq \gamma^2 C \mc{S}^{-1} \int_{\R^N} \big(1+ u_+^{2^*_s-2} \big) (h(u_+))^2 \, dx \nonumber \\
&\leq& \gamma^2 C \mc{S}^{-1} \int_{\R^N} (h(u_+))^2 \, dx + \gamma^2 C \mc{S}^{-1} \int_{\R^N} u_+^{2^*_s-2} (h(u_+))^2 \, dx. \label{eq_dim_stima_buona}
\end{eqnarray}
Let now $R>0$ to be fixed; splitting the second piece of the right-hand side of \eqref{eq_dim_stima_buona} and by using the H\"older inequality we obtain
\begin{align*}
\int_{\R^N} u_+^{2^*_s-2} (h(u_+))^2 \, dx &= \int_{ u \leq R} u_+^{2^*_s-2} (h( u_+))^2 \, dx + \int_{ u>R} u_+^{2^*_s-2} (h( u_+))^2 \, dx \\
&\leq R^{2^*_s-2}\norm{h( u_+)}_2^2 +\left( \int_{ u>R} u^{2^*_s} \, dx \right)^{\frac{2^*_s -2}{2^*_s}}\norm{h( u_+)}_{2^*_s}^2 .
\end{align*}
Since $u \in L^{2^*_s}(\R^N)$, 
we can find a sufficiently large $R=R(\gamma, m_0, \mc{S}^{-1}) $ such that
$$\left( \int_{ u>R} u^{2^*_s} \, dx \right)^{\frac{2^*_s -2}{2^*_s}} < \frac{1}{2} \frac{1}{\gamma^{2} C \mc{S}^{-1}}.$$
Thus, plugging this information into \eqref{eq_dim_stima_buona}, and absorbing the second piece on the right-hand side into the left-hand side, we obtain by \eqref{eq_dim_h_1}
$$\norm{h( u_+)}_{2^*_s}^2 \leq 2\gamma^2 C \mc{S}^{-1} (1+R^{2^*_s-2})\norm{h( u_+)}_2^2 
\leq 2\gamma^2 C \mc{S}^{-1}(1+R^{2^*_s-2})\norm{ u_+}_{2 \gamma }^{2 \gamma}.$$
Recalled that $h=h_{T, \gamma}$, by \eqref{eq_dim_h_3} and Fatou's Lemma we have
\begin{align*}
\norm{ u_+}_{2^*_s \gamma}^{2 \gamma } &= \left( \int_{\R^N} \liminf_{T \to +\infty} h_{T, \gamma}^{2^*_s}( u_+) \, dx \right)^{\frac{2}{2^*_s}} 
\leq \left( \liminf_{T \to +\infty} \int_{\R^N} h_{T, \gamma}^{2^*_s}( u_+) \, dx \right)^{\frac{2}{2^*_s}} \\
&\leq 2\gamma^2 C \mc{S}^{-1} (1+R^{2^*_s-2})\norm{ u_+}_{2 \gamma }^{2 \gamma}.
\end{align*}
By our choice \eqref{eq_dim_stima_gamma} of $\gamma$ we gain that $ u_+ \in L^{2^*_s \gamma}(\R^N)$, which was the claim. By an iteration argument, with 
$$\gamma_0 := \frac{1}{2} 2^*_s, \quad \gamma_i:= \frac{1}{2}2^*_s \gamma_{i-1}, \quad \gamma_i \to +\infty,$$
 we obtain $u_+\in L^r(\R^N)$ for each $r \in [2, +\infty)$. In order to obtain $u_+ \in L^{\infty}(\R^N)$ we need to be careful on the bound on the $L^r$-norms.

Knowing that $u_+$ lies in every Lebesgue space for $r<\infty$ we can implement a more precise iteration argument, where we drop the dependence of the constant on $R$. 
We exploit again \eqref{eq_dim_stima_buona}. Applying again Fatou's Lemma to \eqref{eq_dim_stima_buona} and using \eqref{eq_dim_h_1} we obtain
\begin{equation}\label{eq_nuovo_bootstrap}
\norm{ u_+}_{2^*_s \gamma}^{2 \gamma } \leq \gamma^2 C \mc{S}^{-1}\int_{\R^N}\big(u_+^{2 \gamma} + u_+^{2^*_s-2+2\gamma}\big)\, dx.
\end{equation}
Focusing on the second term on the right-hand side we obtain, exploiting first the generalized H\"older inequality with
$$\frac{1}{N/s} + \frac{1}{2} + \frac{1}{2^*_s} = 1,$$
possible since $u_+^{2^*_s-2} \in L^{\frac{N}{s}}(\R^N)$ because $(2^*_s-2) \frac{N}{s}= \frac{4N}{N-2s}\geq 2$, 
and the generalized Young's inequality then, we obtain
\begin{eqnarray*}
\lefteqn{\int_{\R^N} u_+^{2^*_s-2+2\gamma}\, dx = \int_{\R^N} u_+^{2^*_s-2} u_+^{\gamma} u_+^{\gamma} \, dx \leq \norm{u_+^{2^*_s-2}}_{\frac{N}{s}} \, \norm{u_+^{\gamma}}_2 \, \norm{u_+^{\gamma}}_{2^*_s}} \\
&\leq& \norm{u_+^{2^*_s-2}}_{\frac{N}{s}} \Big( \frac{1}{2 \eps} \norm{u_+^{\gamma}}_2^2+ \frac{\eps}{2} \norm{u_+^{\gamma}}_{2^*_s}^2 \Big) = \norm{u_+}_{\frac{4N}{N-2s}}^{2^*_s-2} \Big( \frac{1}{2 \eps} \norm{u_+}_{2 \gamma}^{2 \gamma}+ \frac{\eps}{2} \norm{u_+}_{2^*_s\gamma}^{2 \gamma} \Big) .
\end{eqnarray*}
Plugging this into \eqref{eq_nuovo_bootstrap}, set $a:= \norm{u_+}_{\frac{4N}{N-2s}}^{2^*_s-2} $, choosing $\eps = \frac{1}{a \gamma^2 C \mc{S}^{-1}}$ and bringing the $L^{2^*_s \gamma}$-norm on the left hand side, we gain
$$ \norm{ u_+}_{2^*_s \gamma}^{2 \gamma } \leq 2 \gamma^2 C \mc{S}^{-1} \big( 1 + \tfrac{1}{2}a^2\gamma^{2} C \mc{S}^{-1} \big) \norm{u_+}_{2\gamma}^{2\gamma} \leq C' \gamma^{4} \norm{u_+}_{2\gamma}^{2\gamma}$$
for some $\gamma$-independent $C'>0$. Choosing $2\gamma_i := 2^*_s \gamma_{i-1}$ we obtain
$$ \norm{ u_+}_{2^*_s \gamma_i} \leq \big(C' \gamma_i^{4}\big)^{\frac{1}{2\gamma_i}} \norm{u_+}_{2^*_s\gamma_{i-1}}$$
and thus
$$ \norm{ u_+}_{2^*_s \gamma_i} 
\leq \prod_{j=0}^i \big(C' \gamma_j^{4}\big)^{\frac{1}{2\gamma_j}} \norm{u_+}_{2^*_s\gamma_0} 
= e^{\sum_{j=0}^i \frac{\log\big(C' \gamma_j^{4}\big)}{2 \gamma_j}} \norm{u_+}_{2^*_s\gamma_0} 
=e^{\sum_{j=0}^i \frac{\log\big(C' \left(\frac{2^*_s}{2}\right)^{4j}\gamma_0^{4}\big)}{2 \left(\frac{2^*_s}{2}\right)^j \gamma_0}} \norm{u_+}_{2^*_s\gamma_0} 
 $$
 and finally, sending $i\to +\infty$,
 $$ \norm{ u_+}_{\infty}
\leq e^{\sum_{j=0}^{\infty}\frac{\log\big(C' \left(\frac{2^*_s}{2}\right)^{4j}\gamma_0^{4}\big)}{2 \left(\frac{2^*_s}{2}\right)^j \gamma_0}} \norm{u_+}_{2^*_s\gamma_0} 
 $$
where the constant is finite. Thus $u_+ \in L^{\infty}(\R^N)$.

To deal with $u_-$ we consider
$$k(t)\equiv k_{T, \gamma}(t):=h_{T, \gamma}(-t), \quad \tilde{k}(t):= \int_t^0 (k'(r))^2 \, dr = \tilde{h}(-t)$$
and choose $\tilde{k}(u)$ as test function. With the same passages we obtain 
$$\norm{k( u)}_{2^*_s}^2 \leq -\mc{S}^{-1} \int_{\R^N} g(x,u) \tilde{k}( u) \, dx$$
and thus
$$\norm{k(u)}_{2^*_s}^2 \leq \mc{S}^{-1}\int_{\R^N} |g(x,u)| \tilde{k}( u) \, dx \leq C \mc{S}^{-1} \int_{\R^N}\big(|u|+ |u|^{2^*_s-1} \big)\tilde{k}( u) \, dx$$
which implies
$$\norm{k(-u_-)}_{2^*_s}^2 \leq C \mc{S}^{-1} \int_{\R^N}\big(|-u_-| +|-u_-|^{2^*_s-1} \big)\tilde{k}(-u_-) \, dx$$
and hence
$$\norm{h(u_-)}_{2^*_s}^2 \leq C \mc{S}^{-1} \int_{\R^N}\big(|u_-| +|u_-|^{2^*_s-1} \big)\tilde{h}(u_-) \, dx;$$
we then proceed as before to gain $u_- \in L^{\infty}(\R^N)$. This concludes the proof.
\QED

\bigskip

\section*{Acknowledgments}
The first and second authors are supported by PRIN 2017JPCAPN ``Qualitative and quantitative aspects of nonlinear PDEs'' and by INdAM-GNAMPA. 
The third author is supported in part by Grant-in-Aid for Scientific Research (19H00644, 18KK0073, 17H02855, 16K13771) of Japan Society for the Promotion of Science.

\section*{Conflict of interest}
All authors declare no conflicts of interest in this paper.



\end{document}